\documentclass[3p]{elsarticle}
\pdfoutput=1
 
\usepackage{amssymb}
\usepackage{amsmath}
\usepackage[all]{xy}
\usepackage{graphicx} 
\usepackage{stmaryrd} 
\usepackage{caption}
\usepackage{subcaption}
\usepackage{mathrsfs}
\usepackage{eucal}
\usepackage{hyperref}
\usepackage{MnSymbol}
\usepackage{latexsym}
\usepackage{dsfont}
\usepackage{adjustbox}
\usepackage{amsthm}
\usepackage{tikz}
\usepackage{enumitem}
\usepackage{bibentry}
\usepackage{float}
\usepackage[normalem]{ulem}

\newcommand{\abs}[1]{\left\vert #1 \right\vert}
\newcommand{\Cov}{\mathrm{Cov}}

\makeatletter
\newmuskip\pFqmuskip

\newcommand*\pFq[6][8]{
  \begingroup
  \pFqmuskip=#1mu\relax
  \mathchardef\normalcomma=\mathcode`,
  \mathcode`\,=\string"8000
  \begingroup\lccode`\~=`\,
  \lowercase{\endgroup\let~}\pFqcomma
  {}_{#2}F_{#3}{\left[\genfrac..{0pt}{}{#4}{#5};#6\right]}
  \endgroup
}
\newcommand{\pFqcomma}{{\normalcomma}\mskip\pFqmuskip}

\begin{document}

\begin{frontmatter}

\title{Tree Pólya Splitting distributions for multivariate count data}

\author[1,2,3,4,5,7]{Samuel Valiquette\corref{cor1}}\ead{samuel.valiquette@mcgill.ca}
\author[3]{Jean Peyhardi}\ead{jean.peyhardi@umontpellier.fr}
\author[5]{Éric Marchand}\ead{eric.marchand@usherbrooke.ca}
\author[3,4]{Gwladys Toulemonde}\ead{gwladys.toulemonde@umontpellier.fr}
\author[1,2,6]{Frédéric Mortier}\ead{frederic.mortier@cirad.fr}

\cortext[cor1]{Corresponding author.}

\affiliation[1]{organization={UPR Forêts et Sociétés, CIRAD},
city={Montpellier},
postcode={34398},
country={France}}

\affiliation[2]{organization={Forêts et Sociétés, Univ Montpellier, CIRAD},
city={Montpellier},
postcode={34398},
country={France}}

\affiliation[3]{organization={IMAG, CNRS, Université de Montpellier},
city={Montpellier},
postcode={34090},
country={France}}

\affiliation[4]{organization={LEMON, Inria},
city={Montpellier},
postcode={34095},
country={France}}

\affiliation[5]{organization={Département de mathématiques, Université de Sherbrooke},
city={Sherbrooke},
postcode={J1K 2R1},
country={Canada}}

\affiliation[6]{organization={AMAP, Univ. Montpellier, IRD, CNRS, CIRAD, INRAE},
city={Montpellier},
country={France}} 

\affiliation[7]{organization={Department of Mathematics and Statistics, McGill University},
city={Montréal},
postcode={H3A 0B9},
country={Canada}}

\begin{abstract}
In this article, we develop a new class of multivariate distributions adapted for count data, called Tree Pólya Splitting.
This class results from the combination of a univariate distribution and singular multivariate distributions along a fixed partition tree.
Known distributions, including the Dirichlet-multinomial, the generalized Dirichlet-multinomial and the Dirichlet-tree multinomial, are particular cases within this class.
As we will demonstrate, these distributions are flexible, allowing for the modeling of complex dependence structures (positive, negative, or null)  at the observation level.
Specifically, we present the theoretical properties of Tree Pólya Splitting distributions by focusing primarily on marginal distributions, factorial moments, and dependence structures (covariance and correlations).
A dataset of abundance of Trichoptera is used, on one hand, as a benchmark to illustrate the theoretical properties developed in this article, and on the other hand, to demonstrate the interest of these types of models, notably by comparing them to other approaches for fitting multivariate data, such as the Poisson-lognormal model in ecology or singular multivariate distributions used in microbiome.
\end{abstract}


\end{frontmatter}

\newtheorem{theorem}{Theorem}
\newtheorem{proposition}[theorem]{Proposition}
\newtheorem{corollary}[theorem]{Corollary}
\newtheorem{lemma}[theorem]{Lemma}
\newdefinition{definition}{Definition}
\newdefinition{rmk}{Remark}
\newproof{pf}{Proof}
\newproof{pot}{Proof of Theorem \ref{thm2}}

\section*{Introduction}\label{Introduction}

Modeling multivariate count data is crucial in many applied fields. In ecology, jointly modeling species distribution according for environmental factors is of primary importance for predicting the impact of climate changes at the ecosystem scale \citep{ovaskainen10,warton15,bry18}. Similar challenges arise in the microbiome context, where understanding microbial community composition may help in defining individual healthcare strategies \citep{Chen2013, Wang2017, Tang2018}, or in econometric analysis to evaluate the number of transactions between various companies \citep{winkelmann2008econometric}.
Finding an appropriate model remains challenging.
In particular, some datasets may exhibit simultaneously positive or negative correlations between different pairs of variables.
An ideal model should be flexible enough to take into account such a correlation structure, while remaining simple for inference and interpretation. 
Further consideration should be given to marginal distributions, which may be overdispersed due to an excess of zeros and/or extreme values present in the sample.

Given these constraints, several models have been proposed, such as the multivariate generalized Waring distribution \citep{Xekalaki}, the discrete Schur-constant model \citep{CASTANER2015343}, and the negative multinomial.
An essential feature of these examples is their representation.
Indeed, as presented by \cite{JONES201983} and \cite{PEYHARDI2021104677}, these models belong to a large class of distributions where each can be expressed as a composition of a univariate discrete distribution and a singular multivariate distribution.
Precisely, \cite{JONES201983} proposed the \textit{sums and shares} model where $\mathbf{Y} = (Y_1, \dots, Y_J) \in \mathbb{N}^J$ is such that the distribution of $\mathbf{Y}$ given $\sum_{j=1}^J Y_j = n$ is Dirichlet-multinomial, and the distribution of $\sum_{j=1}^J Y_j$ is negative binomial.
More general cases were studied by \cite{PEYHARDI2021104677} using singular distributions with certain properties (e.g. multivariate Pólya distributions introduced by \cite{eggenberger1923statistik}), and an arbitrary discrete univariate distribution on the sum.
Intuitively, their models can be interpreted as the random sharing of a univariate random variable into $J$ categories.
This simple stochastic representation where the sum of $\mathbf{Y} \in \mathbb{N}^J$ is randomly split by a Pólya distribution is called the \textit{Pólya Splitting} distribution.
This class emerges naturally as stationary distributions of a multivariate birth–death process under extended neutral theory \citep{Split_neutral}, possessing tractable univariate and multivariate marginals. Its dispersion is well understood, as is the dependence structure. However, Pólya Splitting are restricted since all pairwise correlations necessarily have the same sign.

Many applications consider only the singular multivariate distribution.
This is particularly true in biology, where RNA-sequences are studied.
In this field, the Dirichlet-multinomial and its many generalizations are widely utilized \citep[e.g.][]{Chen2013}.
The generalized Dirichlet-multinomial, introduced by \cite{Mosimann1969}, is considered by \cite{Tang&Chen2018} with a focus on zero-inflation. 
Another example is the Dirichlet-tree multinomial model proposed by \cite{Dennis1991} and employed by \cite{Wang2017} for gut microorganisms.
These examples also have an interesting representation.  
Indeed, for $\mathbf{Y} \in \mathbb{N}^J$ such that $\sum_{j=1}^J Y_j = n$, each distribution can be interpreted as a stochastic process where the total is recursively split by multiple Dirichlet-multinomial distributions.
Such a process can be represented by a tree structure where each node is distributed conditionally as Dirichlet-multinomial, and the leaves are the marginals $Y_j$ $j= 1, \dots, J$.
In \cite{Wang2017}, the authors justified their application of the Dirichlet-tree multinomial with the phylogenetic tree structure of their data.
A similar justification is provided in \cite{Fabrice} for the analysis of forest tree abundance.
This tree-like structure of the distribution enables various sign of correlation, but is less flexible since $\sum_{j=1}^J Y_j$ is fixed and not random.
In fact, as we will show in this work, such a constraint has a significant impact not only on the correlation structure, but also on its marginals.

Alternatively, the multivariate Poisson-lognormal distribution has been proposed as a flexible solution to model both correlations and margins \citep{AITCHISON}.
This model exhibits a diverse correlation structure due to its underlying multivariate lognormal distribution as a latent variable. 
However, from an application perspective, it is important to note that such a dependence structure represents the latent space rather than the observations.
Specifically, the conditional independencies induced by the multivariate Poisson-lognormal distribution are determined by its precision matrix; however, such relationships may not hold for the observations.
A compelling model would match the flexibility of the Poisson-lognormal, while being simpler and adapted to the discrete data by avoiding the use of a continuous latent distribution to describe their correlations. 

In this article, we propose a new class of multivariate discrete distributions named \textit{Tree Pólya Splitting}, which combines a univariate distribution with a tree singular distribution, where each node is associated with a Pólya split.
As it will be demonstrated, this simple modification of the Pólya Splitting enables overdispersed marginals, a diverse correlation structure and a true data structures of dependence.
This composition of sum and tree splitting also allows for a straightforward inference approach where each component is estimated independently.
Since this new model is a generalization of \cite{PEYHARDI2021104677}, we present various properties of the Tree Pólya Splitting and compare them to those of the Pólya Splitting.
This paper is organized as follows.
Section 1 introduces notations and basic results of Pólya Splitting that are used throughout the paper. 
We also provide new results concerning the dispersion of marginal distributions and bounds for correlations.
Section 2 introduces the class of Tree Pólya Splitting distributions.
First, we define the tree structure and associated distributions.
Following this, properties of marginal distributions, factorial moments, and covariance/correlation are presented.
A detailed study is carried out for each property with the help of a running example.
Finally, in Section 3, we present an application of the Tree Pólya Splitting to the Trichoptera dataset provided by \cite{usseglio1987influence} and compare it to the Poisson-lognormal.
We also briefly explore how the observed data can inform us on the underlying tree structure.
All proofs of propositions and corollaries are presented in this paper are given in the Appendix.

\section{Pólya Splitting distributions}\label{section:Pólya_Split}

This section presents the notations, definitions, and properties of the Pólya Splitting distribution used throughout the paper and serving as building blocks for defining the Tree Pólya Splitting distribution.
Marginal distributions, factorial moments, and the Pearson correlation structure of the Pólya Splitting model will be presented.
We refer to \cite{JONES201983,Peyhardi2022,PEYHARDI2021104677,valiquette_thesis} for further details.

\subsection{Definitions and notations} \label{subsection:Def_Not}
Vectors and scalars will be denoted by bold and plain letters, respectively. 
For a vector $\mathbf{y} = (y_1, \dots, y_J)$, the sum of its components is denoted by $\abs{\mathbf{y}} = \sum_{j=1}^J y_j$.
For $\mathcal{J} \subset \{1, \dots, J\}$ a subset of indices, we denote $\mathbf{y}_{\mathcal{J}}$ and its complement $\mathbf{y}_{-\mathcal{J}}$ as $\mathbf{y}_{\mathcal{J}} = (y_j)_{j \in \mathcal{J}}$ and $\mathbf{y}_{-\mathcal{J}} = (y_j)_{ j \in \{1, \dots, J\} \setminus \mathcal{J}}$.
Also, any binary operation between vectors is taken component-wise.
The \textit{discrete simplex} will be denoted by $\Delta_n := \left\{ \mathbf{y} \in \mathbb{N}^J: \abs{\mathbf{y}} = n\right\}$.
For $\theta \in \mathbb{R}_+$, $c \in \{-1,0,1\}$ and $n \in \mathbb{N}_+$, the function $(\theta)_{(n,c)}$ denotes the \textit{generalized factorial} given by
\begin{equation}
    \label{eqn:gfactorial}
    (\theta)_{(n,c)} = \prod_{j=0}^{n-1} (\theta + j c),
\end{equation}
and $(\theta)_{(0,c)}=1$.
If $c = 0$, then $(\theta)_{(n,0)} = \theta^n$, while $c = -1$ and $c=1$ correspond respectively to the \textit{falling} and \textit{rising} factorial.
These will be denoted by the Pochhammer symbols $(\theta)_{(n)} := (\theta)_{(n,-1)}$ and $(\theta)_n := (\theta)_{(n,1)}$ respectively.
Finally, for any $\boldsymbol{\theta} \in \mathbb{R}_+^J$, $\mathbf{r} \in \mathbb{N}^J$ let us denote $(\boldsymbol{\theta})_\mathbf{r} := \prod_{j=1}^J (\theta_j)_{r_j}$.
The special case where $\mathbf{r} = (n, \dots, n)$ for $n \in \mathbb{N}$ will be denoted by $(\boldsymbol{\theta})_n$.
Both notations will also apply to the falling factorial.
A random variable $\mathbf{Y} \in \Delta_n$ is \textit{Pólya} distributed \citep{eggenberger1923statistik, johnson1997discrete} if its probability mass function (p.m.f.) is given by 
\begin{equation}
\label{eqn:Pólya}
    \mathbb{P}_{\abs{\mathbf{Y}}=n}(\mathbf{Y} = \mathbf{y}) = \frac{n!}{(\abs{\boldsymbol{\theta}})_{(n,c)}} \prod_{j=1}^J \frac{(\theta_{j})_{(y_j,c)}}{y_j!} \mathds{1}_{\Delta_n}(\mathbf{y}),
\end{equation}
for $c \in \{-1,0,1\}$ and parameter $\boldsymbol{\theta} = (\theta_1, \dots, \theta_J)$. 
Such a distribution will be denoted by $\mathcal{P}_{\Delta_n}^{[c]}(\boldsymbol{\theta})$.
The following distributions are retrieved: the hypergeometric distribution $\mathcal{H}_{\Delta_n}(\boldsymbol{\theta})$ ($c=-1$), the multinomial $\mathcal{M}_{\Delta_n}(\boldsymbol{\theta})$ ($c=0$), and the Dirichlet-multinomial $\mathcal{DM}_{\Delta_n}(\boldsymbol{\theta})$ ($c=1$).
In order to have an adequate distribution on $\Delta_n$, the allowable values of $\boldsymbol{\theta}$ are the following: $\boldsymbol{\theta} \in \mathbb{R}_+^J$ for $c \in \{0, 1\}$, and $\boldsymbol{\theta} \in \mathbb{N}_+^J$ such that $\abs{\boldsymbol{\theta}} \geq n$ for $c = -1$.
Additionally, the univariate Pólya, i.e. its non-singular form, will be denoted by $\mathcal{P}^{[c]}_n(\theta, \tau)$ and its p.m.f. is given by
\begin{equation}
\label{eqn:uni_Pólya}
    \mathbb{P}(Y = y) = \binom{n}{y} \frac{(\theta)_{(y,c)}(\tau)_{(n-y,c)}}{(\theta + \tau)_{(n,c)}} \mathds{1}_{\{0,\dots,n\}}(y).
\end{equation}

Combining the distribution \eqref{eqn:Pólya} with the hypothesis that $\abs{\mathbf{Y}} \sim \mathcal{L}(\boldsymbol{\psi})$, a univariate discrete distribution, yields a discrete multivariate distribution with $J$ degrees of freedom.
This distribution is named the \textit{Pólya Splitting distribution} and is defined as follows. 

\begin{definition}{(Pólya Splitting distribution)}
\label{def:Splitting}
    A random vector $\mathbf{Y} = (Y_1, \dots, Y_J) \in \mathbb{N}^J$ follows a Pólya Splitting distribution with parameters $c$, $\boldsymbol{\theta}$, $\boldsymbol{\psi}$, and sum distribution $\mathcal{L}(\boldsymbol{\psi})$ if  $\abs{\mathbf{Y}} \sim \mathcal{L}(\boldsymbol{\psi})$ and $\mathbf{Y} \hspace{1 pt} \mid \hspace{1 pt} \abs{\mathbf{Y}} = n  \sim \mathcal{P}_{\Delta_n}^{[c]}(\boldsymbol{\theta})$.
This composition is denoted by 
$$\mathbf{Y} \sim \mathcal{P}^{[c]}_{\Delta_n}(\boldsymbol{\theta}) \underset{n}{\wedge} \mathcal{L}(\boldsymbol{\psi}).$$
Its p.m.f. is given by 
\begin{equation}
    \label{eqn:Splitting_pmf}
    \mathbb{P}(\mathbf{Y} = \mathbf{y}) = \mathbb{P}(\abs{\mathbf{Y}} = n) \left[\frac{n!}{(\abs{\boldsymbol{\theta}})_{(n,c)}} \prod_{j=1}^J \frac{(\theta_{j})_{(y_j,c)}}{y_j!}\right],
\end{equation}
with $n = \abs{\mathbf{y}}$ and $\mathbb{P}(\abs{\mathbf{Y}} = n)$ the p.m.f. of $\mathcal{L}(\boldsymbol{\psi})$.
\end{definition}

Before proceeding, the case $c=-1$ needs to be carefully analyzed. 
Indeed, the restriction on $\abs{\boldsymbol{\theta}} \geq n$ is stated for $n$ fixed.
However, in a Pólya Splitting model, this value is random.
Therefore, it is required that the support of $\mathcal{L}(\boldsymbol{\psi})$ be finite with upper bound value $m \in \mathbb{N}_+$, in which case the Pólya Splitting distribution for $c = - 1$ is well defined only if $\abs{\boldsymbol{\theta}} \geq m$.

\subsection{Marginal distributions and factorial moments}

The marginals of the Pólya Splitting distribution have a particular representation known as Pólya thinning.
Indeed, it is shown in \cite{PEYHARDI2021104677} that for $\mathbf{Y} \sim \mathcal{P}^{[c]}_{\Delta_n}(\boldsymbol{\theta}) \underset{n}{\wedge} \mathcal{L}(\boldsymbol{\psi})$, the marginal distribution of $Y_j$ is given by
\begin{equation}
\label{eqn:marg}
    Y_j \sim \mathcal{P}^{[c]}_n(\theta_j, \abs{\boldsymbol{\theta}_{-j}}) \underset{n}{\wedge} \mathcal{L}(\boldsymbol{\psi}).
\end{equation}
Note that the univariate distribution $\mathcal{L}(\boldsymbol{\psi})$ in \eqref{eqn:marg} is, in a sense, reduced by the univariate Pólya distribution, hence the term Pólya thinning.
This thinning operator is a generalization of the binomial thinning operator studied in \cite{Rao1965} and also used in the time series model proposed by \cite{mckenzie1985some} and \cite{al1987first}.
It can also be further generalized to the so-called quasi-Pólya thinning operator.
For more details of the latter, see \cite{Peyhardi2022}, and for an overview of the thinning operator, see \cite{weiss2008thinning} and \cite{scotto2015thinning}.
\cite{Peyhardi2022} presents three general families of distributions that are stable under this operator, i.e. distributions $\mathcal{L}$ such that
$$\mathcal{L}(\Tilde{\boldsymbol{\psi}}) = \mathcal{P}_n^{[c]}(\theta, \tau) \underset{n}{\wedge} \mathcal{L}(\boldsymbol{\psi})$$
where $\Tilde{\boldsymbol{\psi}}$ is a parameter update of $\boldsymbol{\psi}$.
One such family consists of \textit{power series distribution}, denoted by $\mathcal{PS}^{[c]}(\theta, \alpha)$, with p.m.f. given by
\begin{equation}
\label{eqn:PowerSerie}
    \mathbf{p}(y) \propto \frac{\alpha^y}{y!} (\theta)_{(y,c)},
\end{equation}
where the values of the parameters $\alpha$, $\theta$, and support of $Y$ depend on $c$.
For each value of $c \in \{-1, 0, 1\}$, 
the corresponding distributions are the binomial, Poisson, and negative binomial respectively.
These distributions are stable under the hypergeometric, binomial, and beta-binomial thinning operator respectively.
See Table \ref{Tab:power_dist} for each distribution represented in terms of $\alpha$, $\theta$ and $c$.

\begin{table}[H]
\centering
\begin{adjustbox}{width=0.75\textwidth}
\small
\begin{tabular}{|c|c|c|c|}
\hline
 Distributions & Parameters & Support & P.m.f.\\ 
\hline
   Binomial $\mathcal{B}_\theta(\alpha)$ & $c=-1$, $\theta \in \mathbb{N}_+$, $\alpha \in \mathbb{R}_+$ & $y \in \{0,...,\theta\}$ & $\binom{\theta}{y} \left(\frac{\alpha}{\alpha+1}\right)^y \left(\frac{1}{\alpha+1}\right)^{\theta-y}$\\
\hline
   Poisson $\mathcal{P}(\alpha\theta)$ & $c=0$, $(\theta, \alpha) \in \mathbb{R}^2_+$ & $y \in \mathbb{N}$  &  $\frac{(\alpha \theta)^n}{n!} e^{-\alpha \theta}$  \\
\hline
   Negative binomial $\mathcal{NB}(\theta, \alpha)$ & $c=1$, $\theta \in \mathbb{R}_+$, $\alpha \in(0,1)$ & $y \in \mathbb{N}$  &   $\frac{(\theta)_{y}}{y!} \alpha^y (1-\alpha)^\theta$  \\
 \hline 
\end{tabular}
\end{adjustbox}
\caption{Cases of power series distributions $\mathcal{PS}^{[c]}(\theta, \alpha)$ for $c = \{-1,0,1\}$ and parameters $\theta$ and $\alpha$.}
\label{Tab:power_dist}
\end{table}

For Pólya Splitting, the power series distributions $\mathcal{PS}^{[c]}(\theta, \alpha)$ are the only ones which allow the marginals to be independent.
Indeed, \cite{Peyhardi2022} shows that for $\mathbf{Y} \sim \mathcal{P}_{\Delta_n}^{[c]}(\boldsymbol{\theta}) \underset{n}{\wedge} \mathcal{L}(\boldsymbol{\psi})$, the $Y_j$'s are independent if and only if $$\mathcal{L}(\boldsymbol{\psi}) = \mathcal{PS}^{[c]}(\abs{\boldsymbol{\theta}}, \alpha).$$
Finally, we define the factorial moments, which are particularly useful for random variables in $\mathbb{N}^J$, as their probabilities can be expressed in terms of them \citep{johnson1997discrete, johnson2005univariate}.
For a univariate random variable, the $r$-th factorial moment of $Y \in \mathbb{N}$ is the expected value of the $r$-th falling factorial, i.e. $\mathrm{E}[(Y)_{(r)}]$.
Similarly for $\mathbf{r} = (r_1, \dots, r_J) \in \mathbb{N}^J$, the factorial moment of $\mathbf{Y} \in \mathbb{N}^J$ is the expectation $\mathrm{E}[(\mathbf{Y})_{(\mathbf{r})} ].$
The factorial moments of the Pólya splitting distribution are as follows.

\begin{proposition}
\label{prop:factorial_moment}
    For $\mathbf{Y}\sim \mathcal{P}_{\Delta_n}^{[c]}(\boldsymbol{\theta}) \underset{n}{\wedge} \mathcal{L}(\boldsymbol{\psi})$ and $\mathbf{r} \in \mathbb{N}^J$, the factorial moments are given by 
\begin{equation}
\label{eqn:factorial_moment}
    \mathrm{E}[(\mathbf{Y})_{(\mathbf{r})} ] = \frac{\mu_{\abs{\mathbf{r}}}}{(\abs{\boldsymbol{\theta}})_{(\abs{\mathbf{r}}, c)}} \prod_{j=1}^J (\theta_j)_{(r_j,c)},
\end{equation}
where $\mu_{\abs{\mathbf{r}}}$ is the $\abs{\mathbf{r}}$-th factorial moment of $\mathcal{L}(\boldsymbol{\psi})$. 
As special cases, we have for every $i \neq j$ 
\begin{align}
    \label{eqn:factorial_moment_cases}
\mathrm{E}[Y_i] &= \mu_1 \frac{\theta_i}{\abs{\boldsymbol{\theta}}} \hspace{0.3cm} \mbox{and} \hspace{0.3cm} \mathrm{E}[Y_i Y_j] = \mu_2 \frac{\theta_i \theta_j}{\abs{\boldsymbol{\theta}}(\abs{\boldsymbol{\theta}}+c)}.
\end{align}
\end{proposition}

\subsection{Covariance and dispersion}
\label{section:covariance}

The covariance between $Y_i$ and $Y_j$ for $i \neq j$ in the Pólya Splitting distribution can be deduced by the moments \eqref{eqn:factorial_moment_cases} and is given by
    \begin{equation}
    \label{eqn:covariance}
    \mathrm{Cov}(Y_i, Y_j) = \frac{\theta_i \theta_j}{\abs{\boldsymbol{\theta}}^2(\abs{\boldsymbol{\theta}}+c)} \left[ (\mu_{2} - \mu_{1}^2) \abs{\boldsymbol{\theta}} - c \mu_{1}^2 \right].
\end{equation}
Notice that when $\mathcal{L}(\boldsymbol{\psi}) = \mathcal{D}_m$, a Dirac distribution at $m$, equation \eqref{eqn:covariance} reduces to the covariance of the Pólya distribution, which is always negative for any $c \in \{-1,0,1\}$. 
Here, we are interested in the signs of the covariances as they relate to the distribution $\mathcal{L}(\boldsymbol{\psi})$.
Clearly, the sign of (\ref{eqn:covariance}) is related to the hyperplane defined by 
$(\mu_{2} - \mu_{1}^2) \abs{\boldsymbol{\theta}} = c \mu_{1}^2$ and separates the parameter values $\boldsymbol{\theta}$ into regions of negative, positive, and null covariance.
Observe as well that the covariances have the same sign for all pairs $(Y_i, Y_j)$.
Additionally, note that $\mu_{2} - \mu_{1}^2 = \mathrm{Var}\left[ \abs{\mathbf{Y}}\right] - \mathrm{E}\left[ \abs{\mathbf{Y}}\right]$. 
This determines what type of dispersion the distribution $\mathcal{L}(\boldsymbol{\psi})$ has.
There are three possible situations, $\mathcal{L}(\boldsymbol{\psi})$ is \textit{underdispersed} if $\mu_{2} - \mu_{1}^2 < 0$, \textit{overdispersed} if $\mu_{2} - \mu_{1}^2 > 0$, or \textit{null dispersed} if $\mu_{2} - \mu_{1}^2 = 0$.
Each type of dispersion and value $c$ of the Pólya yields different outcomes.
For $c = 0$, the sign of covariance is simply determined by the dispersion of $\mathcal{L}(\boldsymbol{\psi})$.
In the case of Dirichlet-multinomial Splitting ($c=1$), then
if $\mathcal{L}(\boldsymbol{\psi})$ is underdispersed or has null dispersion, the sign is always negative.
However, if $\mathcal{L}(\boldsymbol{\psi})$ is overdispersed, the sign of covariance is negative, null or positive if and only if $\abs{\boldsymbol{\theta}}$ is less, equal or greater than ${\mu_{1}^2}/{(\mu_{2}-\mu_{1}^2)}$ respectively.
A similar analysis for $c= - 1$ can be made using the restriction on $\boldsymbol{\theta}$.

Since the dispersion of the distribution $\mathcal{L}(\boldsymbol{\psi})$ is relevant for the covariance sign, it will be useful for the new model to understand how this dispersion is preserved in the marginals.
For example, if $\mathcal{L}(\boldsymbol{\psi})$ is overdispersed, does it imply that the marginals are necessary overdispersed?
We have the following.

\begin{proposition}
\label{prop:dispersion}
    For $\mathbf{Y}\sim \mathcal{P}_{\Delta_n}^{[c]}(\boldsymbol{\theta}) \underset{n}{\wedge} \mathcal{L}(\boldsymbol{\psi})$, then:
    \begin{itemize}
        \item If $c = 0$, the marginals have the same type of dispersion as $\mathcal{L}(\boldsymbol{\psi})$;
        \item If $c = 1$ and $\mathcal{L}(\boldsymbol{\psi})$ has null or positive dispersion, then the marginals are overdispersed;
        \item If $c = -1$ and $\mathcal{L}(\boldsymbol{\psi})$  has null or negative dispersion, then the marginals are underdispersed.
    \end{itemize}
\end{proposition}

Furthermore, Proposition \ref{prop:dispersion} implies that the type of dispersion is preserved for $c = 0$, but can change otherwise. 
For instance, if $\mathcal{L}(\boldsymbol{\psi})$ is underdispersed and $c=1$, then it is possible to have different type of dispersion for the marginals.

\subsection{Pearson correlation structure}

An interesting way to formulate the correlation between $Y_i$ and $Y_j$ is to use the relation between factorial moments of $\mathcal{L}(\boldsymbol{\psi})$ and the marginals.
Indeed, by Proposition \ref{prop:factorial_moment}, it is possible to express $\mathrm{Var}[Y_i]$ in terms of $\Cov(Y_i,Y_j)$ and the moments of $Y_i$, similarly for $\mathrm{Var}[Y_j]$.
The following proposition shows such a relation.

\begin{proposition}
    \label{prop:correlation}
    For $\mathbf{Y}\sim \mathcal{P}_{\Delta_n}^{[c]}(\boldsymbol{\theta}) \underset{n}{\wedge} \mathcal{L}(\boldsymbol{\psi})$ and $i \neq j$, then the Pearson correlation coefficient is given by
    \begin{equation}
        \label{eqn:corr}
        \mathrm{Corr}(Y_i, Y_j) = \mathrm{sgn}(1-M_i)\sqrt{\frac{\theta_i \theta_j}{(\theta_i + c)(\theta_j + c)} (1 - M_i) (1 - M_j)} ,
    \end{equation}
    
    where, for $\mu_{r}$ the $r$-th factorial moment of $\mathcal{L}(\boldsymbol{\psi})$, $$M_k = \frac{\mathrm{E}[Y_k]}{\mathrm{Var}[Y_k]} \left( 1 + \frac{c}{\theta_k} \mathrm{E}[Y_k] \right) = \frac{\mu_{1} \left( 1 + \frac{c}{\abs{\boldsymbol{\theta}}} \mu_{1} \right)}{\mu_{2} \left( \frac{\theta_k + c}{\abs{\boldsymbol{\theta}} + c}\right) + \mu_{1} \left( 1 - \mu_{1} \frac{\theta_k}{\abs{\boldsymbol{\theta}}} \right)}, \text{ } k=i,j,$$
    and $\mathrm{sgn}(x) = \mathds{1}_{[0,\infty)}-\mathds{1}_{(-\infty,0]}$. 
\end{proposition}

It was shown in \cite{JONES201983} that for any distribution $\mathcal{L}(\boldsymbol{\psi})$, $\mathrm{Corr}(Y_i, Y_j) < 1/2$ when $\boldsymbol{\theta} = (1,\dots,1)$ and $c = 1$.
With Proposition \ref{prop:correlation}, we can generalize their result to other values of $\boldsymbol{\theta}$.

\begin{proposition}
    \label{prop:bound_corr}
    For $\mathbf{Y}\sim \mathcal{DM}_{\Delta_n}(\boldsymbol{\theta}) \underset{n}{\wedge} \mathcal{L}(\boldsymbol{\psi})$ and $i \neq j$, then the correlation is such that $$\mathrm{Corr}(Y_i, Y_j) < \sqrt{\frac{\theta_i \theta_j}{(\theta_i + 1)(\theta_j + 1)}} \hspace{5pt}.$$
\end{proposition}

Notice that this bound is not sharp. 
Interestingly, this bound is equal to the geometric mean of $\theta_i/(\theta_i+1)$ and $\theta_j/(\theta_j+1)$.

\section{Tree Pólya Splitting Distribution}

In this section, we first present notations and definitions for rooted trees inspired by \cite{Tang2018}.
Following this, we define the Tree Pólya Splitting distribution and present properties which parallel and generalize those of the previous section.
A running example is used throughout to illustrate and further explore these results.
In particular, we are able to obtain a marginal p.m.f. that generalizes the one given by \cite{JONES201983}.
We also show how the pairwise correlations of the Tree Pólya distribution can indeed take various signs.

\subsection{Definitions and Notations}
\label{subsection:Def_Tree}

Let $\mathfrak{T} =$ $(\mathcal{N},\mathcal{E})$ be an undirected graph with nodes $\mathcal{N}$ and edges $\mathcal{E}$.
$\mathfrak{T}$ is a \textit{tree} if it is connected, i.e. there is a chain of edges between every pair of nodes in the graph, and acyclic, i.e. the graph contains no cycle.
Furthermore, $\mathfrak{T}$ is a \textit{rooted tree} if it is a tree with a fixed node named root, denoted by $\Omega$.
Fixing such a root induces an orientation on every edge in $\mathfrak{T}$.
Precisely, we can recursively establish a parent/child relation between the two nodes that constitute each edge.
For any node $A \in \mathcal{N}$, its set of \textit{children nodes} contains any node directly connected to $A$ in the opposite direction of $\Omega$.
Such a set is denoted by $\mathfrak{C}_A = \{C_1, \dots, C_{J_A}\}$ with $J_A$ the number of children. 
Similarly, the parent of $A$ is the node directly connected in the direction of $\Omega$.
It is denoted by $\mathscr{P}(A)$ and is such that: (i) $\mathscr{P}(\Omega) = \emptyset$, and (ii) $\mathscr{P}(C_i) = \mathscr{P}(C_j) = A$ for $C_i, C_j \in \mathfrak{C}_A$ with $i \neq j$, i.e. $C_i$ and $C_j$ are \textit{sibling nodes}.
Finally, a \textit{leaf} is a node such that its children set is empty, and an \textit{internal node} is any node that is not a leaf. 
In this instance, $\mathfrak{L} \subseteq \mathcal{N}$ and $\mathfrak{I} \subseteq \mathcal{N}$ will denote the subsets of leaves and internal nodes respectively.
Based on these definitions and notations, we are now able to define a specific type of rooted tree useful for the model.

\begin{definition}[Partition tree]
\label{def:Partition_tree}
A rooted tree $\mathfrak{T}$ is a partition tree if its root is given by $\Omega = \{1, \dots, J\}$,
the leaves by $\mathfrak{L} = \{\{1\},\dots,\{J\}\}$, and each sibling forms a non-trivial partition of their parent.
\end{definition}

Notice that Definition \ref{def:Partition_tree} implies that any internal node $A \in \mathfrak{I}$ in a partition tree must have at least two children.
To understand various properties related to the Tree Pólya distribution, we need to study parts of the partition tree $\mathfrak{T}$.
First, we need to define subtrees of $\mathfrak{T}$, named \textit{pruned trees}.

\begin{definition}[Pruned Tree at $A$]
\label{def:pruned_tree}
For any partition tree $\mathfrak{T}$ and $A \in \mathfrak{I}$, $\mathfrak{T}_A$ is the pruned tree of $\mathfrak{T}$ at $A$ if its set of internal nodes and leaves are given by $\mathfrak{I}_A = ( B \in \mathfrak{I}: B \subseteq A)$ and $\mathfrak{L}_A = (\{j\} \in \mathfrak{L}: \{j\} \subseteq  A)$ respectively. 
\end{definition}

Finally, the notion of path between two nodes is needed.
Such a path is constructed through an iteration of the parent nodes from any leaf or internal node to another node.

\begin{definition}[Path]
\label{def:path}
For any partition tree $\mathfrak{T}$ and $A \in \mathfrak{I} \cup \mathfrak{L}$, let
$$\mathscr{P}^n_{A} := \underbrace{(\mathscr{P} \circ \mathscr{P} \circ \cdots \circ \mathscr{P})}_{\text{n times}}(A)$$
be the $n$-th parent of node $A$ with $\mathscr{P}^0_{A} = A$. 
For $B \in \mathfrak{I}$ such that $A \subset B$, the path from $A$ to $B$ is defined by the ordered set $$ \mathrm{Path}^B_A := \left(\mathscr{P}^0_{A}, \mathscr{P}^1_{A}, \dots, \mathscr{P}^K_{A}\right),$$ where $K$ is such that $\mathscr{P}_A^K = B$.
By convention, $A_n \in \mathrm{Path}_A^B$ means that $A_n$ is the $n$-th element of $\mathrm{Path}_A^B$.
Therefore, the element $A_{n-1}$ for $n \geq 1$ should be interpreted as the child of $A_n$.
Moreover, if $B = \Omega$, then we note $\mathrm{Path}_A:=\mathrm{Path}_A^\Omega$.
\end{definition}

With these definitions, the structure of the partition tree can be fully described and used. 
For the running example, we will use the partition tree presented in Figure \ref{fig:example}.
In this example, $\Omega = \{1,\dots,10\}$, $A =\{4,\dots, 10\}$ is an internal node with children nodes $\mathfrak{C}_A = \{\{4,5\}, \{6,7\}, \{8,9,10\}  \}$, $\mathscr{P}(A) = \Omega$, and the path between the leaf $\{9\}$ and $A$ is given by $\mathrm{Path}_{\{4\}}^A = (\{9\}, \{9,10\}, \{8,9,10\}, A)$.

\begin{figure}[ht]
\centering
\begin{subfigure}{.5\textwidth}
  \centering
    \resizebox{0.9\textwidth}{0.5\textwidth}{
\begin{tikzpicture}[scale = 0.9]
    \node {$\{1,\dots,10\}$}[sibling distance = 3.5cm]
            child {node {$\{1,2\}$}[sibling distance = 2cm]
                child{node{$\{1\}$}}
                child{node{$\{2\}$}}
             }
             child {node {$\{3\}$}
            }
            child {node {$\{4,\dots,10\}$}[sibling distance = 3cm]
      child {node {$\{4,5\}$}[sibling distance = 1.5cm]
                child{node{$\{4\}$}}
                child{node{$\{5\}$}}
      }
      child {node {$\{6,7\}$}[sibling distance = 1.5cm]
               child{node{$\{6\}$}}
                child{node{$\{7\}$}}
      }
      child {node {$\{8,9,10\}$}[sibling distance = 1.5cm]
                child{node{$\{8\}$}}
                child{node{$\{9,10\}$}
                    child{node{$\{9\}$}}
                    child{node{$\{10\}$}}            
                }
      }
    };
\end{tikzpicture}
}
  \caption{Partition tree with 10 leaves and 7 internal nodes.}
  \label{fig:example}
\end{subfigure}%
\begin{subfigure}{.5\textwidth}
\centering
    \resizebox{0.9\textwidth}{0.6\textwidth}{
\begin{tikzpicture}[scale = 0.9]
    \node{$\mathcal{NB}(\alpha,p)$}
    child{ node {$\mathcal{M}$}[sibling distance = 3.5cm]
            child {node[label={[xshift = 1.4cm,yshift= 0.4cm] $\pi_{1,2}$}] {$\mathcal{DM}$}[sibling distance = 2cm]
                child{node[label={[xshift = 0.1cm,yshift= 0.2cm] $\theta_{1}$}]{$Y_1$}}
                child{node[label={[xshift = -0.1cm,yshift= 0.2cm] $\theta_{2}$}]{$Y_2$}}
            }
            child {node[label={[xshift = -0.4cm,yshift= 0.1cm] $\pi_{3}$}] {$Y_3$}
            }
            child {node[label={[xshift = -1.4cm,yshift= 0.4cm] $\pi_{4:10}$}] {$\mathcal{DM}$}[sibling distance = 3cm]
      child {node[label={[xshift = 0.8cm,yshift= 0.2cm] $\theta_{4,5}$}] {$\mathcal{M}$}[sibling distance = 1.5cm]
                child{node[label={[xshift = -0.1cm,yshift= 0.2cm] $\pi_{4}$}]{$Y_4$}}
                child{node[label={[xshift = 0.1cm,yshift= 0.2cm] $\pi_{5}$}]{$Y_5$}}
      }
      child {node[label={[xshift = -0.45cm,yshift= 0cm] $\theta_{6,7}$}] {$\mathcal{DM}$}[sibling distance = 1.5cm]
                child{node[label={[xshift = -0.1cm,yshift= 0.2cm] $\theta_{6}$}]{$Y_6$}}
                child{node[label={[xshift = 0.1cm,yshift= 0.2cm] $\theta_{7}$}]{$Y_7$}}
      }
      child {node[label={[xshift = -0.7cm,yshift= 0.2cm] $\theta_{8,9,10}$}] {$\mathcal{DM}$}[sibling distance = 1.5cm]
                child{node[label={[xshift = -0.1cm,yshift= 0.2cm] $\theta_{8}$}]{$Y_8$}}
                child{node[label={[xshift = 0.1cm,yshift= 0.2cm] $\theta_{9,10}$}]{$\mathcal{M}$}
                    child{node[label={[xshift = 0cm,yshift= 0.2cm] $\pi_{9}$}]{$Y_9$}}
                    child{node[label={[xshift = 0cm,yshift= 0.2cm] $\pi_{10}$}]{$Y_{10}$}}            
                }
      }
    }};
\end{tikzpicture}
}
  \caption{Tree Pólya Splitting distribution using partition tree \ref{fig:example}, $\mathcal{L}(\boldsymbol{\psi}) = \mathcal{NB}(\alpha, p)$, and multinomial/Dirichlet-multinomial splits.}
  \label{fig:example_3}
\end{subfigure}
\caption{Representations of the running example}
\label{fig:ancestral_def}
\end{figure}

It is now possible to generalize the Pólya Splitting distribution with the partition tree as follows.
For $\mathbf{Y} \in \mathbb{N}^J$, $\abs{\mathbf{Y}}$ is split by a Pólya into subsums which are again split until each marginal $Y_j$ is obtained.
It is assumed that these divisions are fixed by the model, i.e. we know the structure of $\mathcal{T}$.
This approach allows to create various clusters of $\mathbf{Y}$ with different dependence structures or marginals.
Since the order of divisions is fixed, this new distribution can be represented and studied with the partition tree $\mathfrak{T}$.
Indeed, the internal nodes $\mathfrak{I}$ determine all the subsums involved and the leaves $\mathfrak{L}$ represent all marginals.
Moreover, by construction, the p.m.f. can be directly obtained as the product of probabilities given each subsums.
We have the following definition.

\begin{definition}[Tree Pólya Splitting distribution]
\label{def:Tree_Pólya}
For a partition tree $\mathfrak{T}$, $\mathbf{Y} \in \mathbb{N}^J$ is said to follow a \textit{Tree Pólya Splitting} distribution if, for each node $A \in \mathfrak{I}$, the distribution of the subsums $\big(|\mathbf{Y}_{C_1}|,\dots, |\mathbf{Y}_{C_{J_A}}|\big)$ given $\abs{\mathbf{Y_A}} = n$ is $\mathcal{P}_{\Delta_n}^{[c_A]}\left( \boldsymbol{\theta}_A \right)$, with $\boldsymbol{\theta}_A = \{\theta_{C}\}_{C \in \mathfrak{C}_A}$ and $c_A$ depending on $A$.
Such a distribution is denoted by $\mathbf{Y} \sim \mathcal{TP}_{\Delta_n}(\mathfrak{T}; \boldsymbol{\theta}, \mathbf{c}) \underset{n}{\wedge} \mathcal{L}(\boldsymbol{\psi})$ with $\boldsymbol{\theta} = \{ \boldsymbol{\theta}_A \}_{A\in\mathfrak{I}}$, $\mathbf{c} = \{ c_A \}_{A\in\mathfrak{I}}$ and p.m.f. 
\begin{equation}
    \label{eqn:tree_polya}
    \mathbb{P}(\mathbf{Y}=\mathbf{y}) =  \mathbb{P}(\abs{\mathbf{Y}} = n) \prod_{A \in \mathfrak{I}} \frac{n_A!}{(\abs{\boldsymbol{\theta}_A})_{(n_A,c_A)}} \prod_{C \in \mathfrak{C}_A}\frac{(\theta_C)_{(n_C, c_A)}}{n_C!},
\end{equation}
where $n_A :=\abs{\mathbf{y}_A}$ for any node $A$, and $n_\Omega := n =|\boldsymbol{y}|$ by definition.
\end{definition}

This recursive definition of the p.m.f. based on a partition tree is similar to those given by \cite{peyhardi2016partitioned} for a categorical response, that corresponds to the simple case $n=1$.  
Otherwise, notice that if the root is the only one internal node ($\mathfrak{I} = \Omega$) in Definition \ref{def:Tree_Pólya}, then the basic Pólya Splitting p.m.f. \eqref{eqn:Splitting_pmf} is recovered in \eqref{eqn:tree_polya}.
Moreover, the largest number of parameters needed can be attained when $\mathcal{T}$ has a binary tree structure, i.e., each internal node has two children.
Therefore, depending on the type of Splittings, $\mathcal{TP}_{\Delta_n}(\mathfrak{T}; \boldsymbol{\theta}, \mathbf{c})$ number of parameters varies between $\abs{\boldsymbol{\psi}} + (J-1)$ and $\abs{\boldsymbol{\psi}} + 2(J-1)$.
Now, using the partition tree in Figure \ref{fig:example}, Figure \ref{fig:example_3} presents the model representation of the running example where all internal nodes are either a multinomial ($\mathcal{M}$) or a Dirichlet-multinomial ($\mathcal{DM}$), and each edge is associated to a parameter of the given Pólya.
Moreover, the distribution of $\abs{\mathbf{Y}}$ is indicated at the top of the tree.
In this case, $\mathcal{L}(\boldsymbol{\psi}) = \mathcal{NB}(\alpha,p)$, the negative binomial.
The distribution of $\mathbf{Y}$ is thus fully specified in Figure \ref{fig:example_3}.

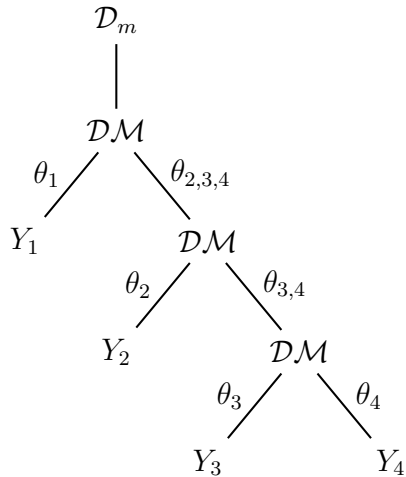
\begin{figure}[ht]
    \centering
        \resizebox{0.3\textwidth}{0.25\textheight}{
\begin{tikzpicture}[thick]
\node{$\mathcal{D}_m$}
child{
    node {$\mathcal{DM}$} [sibling distance = 2.5cm]
    child{node[label={[xshift = 0.3cm,yshift= 0.25cm] $\theta_1$}] {$Y_1$}}
    child{
        node[label={[xshift = -0.1cm,yshift= 0.25cm] $\theta_{2,3,4}$}] {$\mathcal{DM}$} 
        child{node[label={[xshift = 0.3cm,yshift= 0.25cm] $\theta_2$}] {$Y_2$}}
        child{
            node[label={[xshift = -0.2cm,yshift= 0.25cm] $\theta_{3,4}$}] {$\mathcal{DM}$}
            child{node[label={[xshift = 0.3cm,yshift= 0.25cm] $\theta_3$}] {$Y_3$}}
            child{node[label={[xshift = -0.3cm,yshift= 0.25cm] $\theta_4$}] {$Y_4$}}
            }
        }
};
\end{tikzpicture} 
}
    \caption{Generalized Dirichlet-multinomial distribution represented by a Tree Pólya Splitting model with $\mathcal{L}(\boldsymbol{\psi}) = \mathcal{D}_m$, the Dirac distribution at $m$, and a binary cascade tree structure.}
    \label{fig:example_4}
\end{figure}

Just like Pólya Splitting, several known distributions are particular cases of Tree Pólya Splitting.
As previously indicated, the Pólya Splitting itself is a trivial case.
For a fixed value of $\abs{\mathbf{Y}}$, i.e. the total follows a Dirac, the generalized Dirichlet-multinomial can be directly retrieved \citep{Mosimann1969}.  
Indeed, the tree $\mathfrak{T}$ in this case is such that the elements of $\mathfrak{I}$ are given by $A_j = \{j,\dots,J\}$ for all $j\in \{1,\dots,J-1\}$ and their set of children is given by $\mathfrak{C}_{A_j} = \{\{j\}, \{j+1, \dots, J\}\}$.
For each node, a Dirichlet-multinomial is used and can be represented by a binary cascade tree (Figure \ref{fig:example_4}).
Similarly, the Dirichlet-tree multinomial distribution \citep{Dennis1991} uses a more general tree structure where each internal node is again distributed as Dirichlet-multinomial.
In all of the above examples, it is important to keep in mind that the value $\abs{\mathbf{Y}}$ is fixed and not random.
As one may anticipate, the added randomness of $\abs{\mathbf{Y}}$ can have a significant impact on the covariance structure.
Additionally, the order of the marginals is important for any Tree Pólya Splitting.
Indeed, for $\mathbf{Y} = (Y_1, \dots, Y_J) $ and $\Tilde{\mathbf{Y}} = (\mathbf{Y}_{\sigma(1)}, \dots, \mathbf{Y}_{\sigma(J)}) $ with $\sigma(\cdot)$ a non-identity permutation, the tree Splitting p.m.f. is such that $\mathbb{P}(\mathbf{Y} = \mathbf{y}) \neq \mathbb{P}(\Tilde{\mathbf{Y}} = \Tilde{\mathbf{y}})$.
Therefore, the order of the leaves $\mathfrak{L}$ in the tree should always be kept in mind.
This consideration is often overlooked, particularly when modeling species with the generalized Dirichlet-multinomial, where no specific order is defined a priori.

\subsection{Properties}

We expand here on properties of Tree Pólya Splitting distributions, which extend those of Pólya Splitting.
To illustrate them, we use the running example presented in Figure \ref{fig:example_3}, where the total follows a negative binomial distribution, and each internal node can either be a multinomial or a Dirichlet-multinomial distribution.
To understand the impact of $\mathcal{L}(\boldsymbol{\psi})$, we compare the example to the same Tree Pólya Splitting but with a fixed total.

\subsubsection{Marginals}

Since the Tree Pólya Splitting distribution is simply an iteration of different splittings throughout the partition tree, the marginals should have similar form to \eqref{eqn:marg}.
Because any marginal $Y_j$ is represented by a leaf in the partition tree, its path to the root must dictate the form of its distribution as shown in the following proposition.

\begin{proposition}[Univariate marginal for a leaf]
\label{prop:marg_tree_univ}
For $\mathbf{Y} \sim \mathcal{TP}_{\Delta_n}(\mathfrak{T}; \boldsymbol{\theta}, \mathbf{c}) \underset{n}{\wedge} \mathcal{L}(\boldsymbol{\psi})$, the distribution of $Y_j$ is given by
\begin{equation}
    \label{eqn:tree_marg}
    Y_j \sim \bigwedge_{k = 1}^K \mathcal{P}^{\left[c_{A_k}\right]}_{n_k}\left( \theta_{A_{k-1}}, \abs{\boldsymbol{\theta}_{A_k\setminus A_{k-1}}} \right)  \underset{n_{K}}{\wedge} \mathcal{L}(\boldsymbol{\psi}),
\end{equation}
where $A_k \in \mathrm{Path}_{\{j\}}$, $K$ is the length of the path, $\boldsymbol{\theta}_{A_k\setminus A_{k-1}}$ is the set of parameters at node $A_k$ with $\theta_{A_{k-1}}$ removed, and $\bigwedge_{k = 1}^K \mathcal{P}^{\left[c_{A_k}\right]}_{n_k}$ is the iterated composition given by
\begin{equation*}
    \mathcal{P}^{\left[c_{A_1}\right]}_{n_1}\left( \theta_{A_{0}}, \abs{\boldsymbol{\theta}_{A_1\setminus A_{0}}}\right) \underset{n_1}{\wedge} \dots \underset{n_{K-1}}{\wedge}  \mathcal{P}^{\left[c_{A_K}\right]}_{n_K}\left( \theta_{A_{K-1}}, \abs{\boldsymbol{\theta}_{A_K\setminus A_{K-1}}}\right).
\end{equation*}
\end{proposition}

Notice that any partial sum must have a similar distribution since it is represented by an internal node that has a path to the root.
Therefore, for any $A \in \mathfrak{I}$, we can easily adapt Proposition \ref{prop:marg_tree_univ} for $\abs{\mathbf{Y}_{A}}$ by changing $A_k \in \mathrm{Path}_{A}$ in \eqref{eqn:tree_marg}.
Based on this observation, we can derive multivariate marginal distributions that remain consistent with the whole tree structure. 

\begin{proposition}[Multivariate marginal of a subtree]
\label{prop:marg_tree_multi}
For $\mathbf{Y} \sim \mathcal{TP}_{\Delta_n}(\mathfrak{T}; \boldsymbol{\theta}, \mathbf{c}) \underset{n}{\wedge} \mathcal{L}(\boldsymbol{\psi})$ and an internal node $A \in \mathfrak{I}$, the multivariate marginal distribution $\mathbf{Y}_A$ is again Tree Pólya Splitting, i.e. $$\mathbf{Y}_A \sim \mathcal{TP}_{\Delta_n}(\mathfrak{T}_A; \Tilde{\boldsymbol{\theta}},
\Tilde{\mathbf{c}}) \underset{n}{\wedge} \abs{\mathbf{Y}_{A}},$$ where the distribution $\abs{\mathbf{Y}_{A}}$ is given by Proposition \ref{prop:marg_tree_univ}, and $\Tilde{\boldsymbol{\theta}}$, $\Tilde{\mathbf{c}}$ are the parameters involved in $\mathfrak{T}_A$, the pruned tree at $A$.
\end{proposition}

\noindent \textbf{Running example}

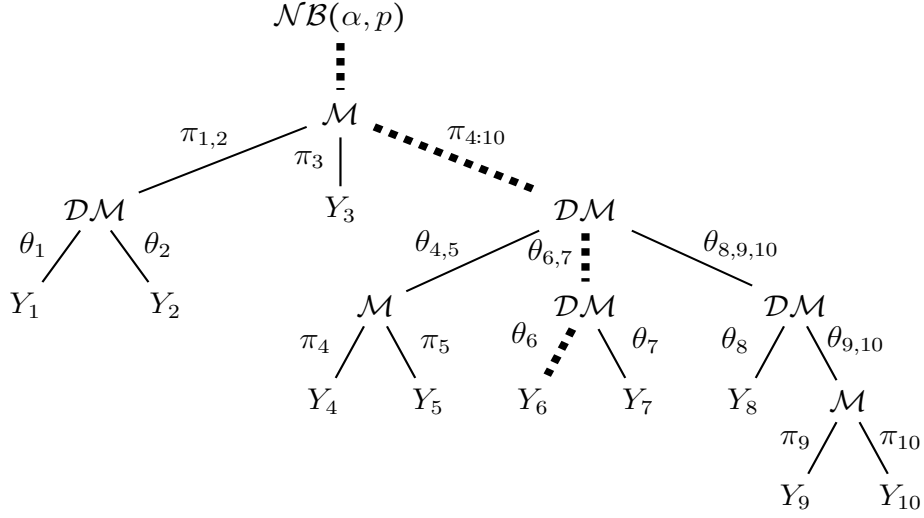
\begin{figure}[H]
    \centering
    \resizebox{0.5\textwidth}{0.2\textheight}{
\begin{tikzpicture}[scale = 0.9]
    \node{$\mathcal{NB}(\alpha,p)$}
    child[line width= 3 pt, dashed]{ node {$\mathcal{M}$}[sibling distance = 3.5cm]
            child[thick, solid] {node[label={[xshift = 1.4cm,yshift= 0.4cm] $\pi_{1,2}$}] {$\mathcal{DM}$}[sibling distance = 2cm]
                child{node[label={[xshift = 0.1cm,yshift= 0.2cm] $\theta_{1}$}]{$Y_1$}}
                child{node[label={[xshift = -0.1cm,yshift= 0.2cm] $\theta_{2}$}]{$Y_2$}}
            }
            child[thick, solid] {node[label={[xshift = -0.4cm,yshift= 0.1cm] $\pi_{3}$}] {$Y_3$}
            }
            child {node[label={[xshift = -1.4cm,yshift= 0.4cm] $\pi_{4:10}$}] {$\mathcal{DM}$}[sibling distance = 3cm]
      child[thick, solid] {node[label={[xshift = 0.8cm,yshift= 0.2cm] $\theta_{4,5}$}] {$\mathcal{M}$}[sibling distance = 1.5cm]
                child{node[label={[xshift = -0.1cm,yshift= 0.2cm] $\pi_{4}$}]{$Y_4$}}
                child{node[label={[xshift = 0.1cm,yshift= 0.2cm] $\pi_{5}$}]{$Y_5$}}
      }
      child {node[label={[xshift = -0.45cm,yshift= 0cm] $\theta_{6,7}$}] {$\mathcal{DM}$}[sibling distance = 1.5cm]
                child{node[label={[xshift = -0.1cm,yshift= 0.2cm] $\theta_{6}$}]{$Y_6$}}
                child[thick, solid]{node[label={[xshift = 0.1cm,yshift= 0.2cm] $\theta_{7}$}]{$Y_7$}}
      }
      child[thick, solid] {node[label={[xshift = -0.7cm,yshift= 0.2cm] $\theta_{8,9,10}$}] {$\mathcal{DM}$}[sibling distance = 1.5cm]
                child{node[label={[xshift = -0.1cm,yshift= 0.2cm] $\theta_{8}$}]{$Y_8$}}
                child{node[label={[xshift = 0.1cm,yshift= 0.2cm] $\theta_{9,10}$}]{$\mathcal{M}$}
                    child{node[label={[xshift = 0cm,yshift= 0.2cm] $\pi_{9}$}]{$Y_9$}}
                    child{node[label={[xshift = 0cm,yshift= 0.2cm] $\pi_{10}$}]{$Y_{10}$}}            
                }
      }
    }};
\end{tikzpicture}
}
    \caption{Path representation of the marginal $Y_6$ in the running example.}
    \label{fig:example_6}
\end{figure}

\noindent As an example, consider the distribution of $Y_6$ in Figure \ref{fig:example_3}.
Using Proposition \ref{prop:marg_tree_univ}, the marginal distribution is given by
\begin{equation}
    \label{eqn:marginal}
    Y_6 \sim \mathcal{BB}_{n_1}(\theta_{6}, \theta_{7}) \underset{n_1}{\wedge} \mathcal{BB}_{n_2}(\theta_{6,7}, \theta_{4,5} + \theta_{8,9,10}) \underset{n_2}{\wedge} \mathcal{B}_{n_3}(\pi_{4:10}) \underset{n_3}{\wedge} \mathcal{NB}(\alpha,p)
\end{equation}
as represented by the path in Figure \ref{fig:example_6}.

\noindent In fact, equation \eqref{eqn:marginal} can be expressed with a composition of only beta-binomial distributions.
Indeed, if the marginal distribution is the composition of binomial, beta-binomial and negative binomial distributions, then all the binomial distributions can be "absorbed" in the negative binomial since the latter is stable under the binomial thinning operator \citep{Peyhardi2022}.
The following formalizes this result.

\pagebreak

\begin{proposition}
    \label{prop:absorb}
    Suppose for $K$ Pólya distributions there are $M$ cases with $c_k = 1$ and parameters $a_k, b_k \in \mathbb{R}_+$, and $K - M$ cases with $c_k = 0$ and parameters $\pi_k \in (0,1)$. 
    Then
    
    $$\bigwedge_{k = 1}^K \mathcal{P}^{[c_k]}_{n_k}\left( \theta_k, \tau_k \right) \underset{n_K}{\wedge} \mathcal{NB}(\alpha,p) = \left[\bigwedge_{m = 1}^M \mathcal{BB}_{n_m}\left( a_m, b_m \right)\right] \underset{n_M}{\wedge} \mathcal{NB}\left(\alpha, \frac{p \gamma}{1-p(1-\gamma)}\right),$$ 
    where $\gamma = \prod_{k=1}^{K-M} \pi_k$, and $(\theta_k, \tau_k)$ is given by $(\pi_k, 1-\pi_k)$ or $(a_k, b_k)$, depending whether $c_k = 0$ or $c_k=1$ respectively.
\end{proposition}

Therefore, to obtain the p.m.f. of \eqref{eqn:marginal}, or any marginal in Figure \ref{fig:example_6}, it is sufficient to study the p.m.f. of the general composition 

\begin{equation}
    \label{eqn:marginal_nb}
    X \sim \left[ \bigwedge_{k = 1}^K \mathcal{BB}_{n_k}\left( a_k, b_k  \right) \right] \underset{n_{K}}{\wedge} \mathcal{NB}(\alpha,p).
\end{equation}
For $K=1$, \cite{JONES201983} and \cite[supplementary material]{PEYHARDI2021104677} showed that the p.m.f. of \eqref{eqn:marginal_nb} is given by 
\begin{equation*}
    (1-p)^\alpha \frac{(\alpha)_n(\alpha_1)_n}{(\alpha_1 + \beta_1)_n} \frac{p^n}{n!} \pFq{2}{1}{\alpha+n, \beta_1}{\alpha_1+\beta_1+n}{p}; \text{ $n \in \mathbb{N}$,}
\end{equation*}
where $\pFq{2}{1}{a,b}{c}{z}$ is \textit{Gauss' hypergeometric function}. 
This result can be generalized for any positive integer $K$.

\begin{proposition}
\label{prop:marginal_nb}
    Let $p \in (0,1)$, $\alpha >0$, $\boldsymbol{a} = (a_1, \dots, a_K)$ and $\boldsymbol{b} = (b_1,\dots,b_K)$ two positive vectors with $K \geq 1$.
    If $X$ is distributed as in \eqref{eqn:marginal_nb}, then its p.m.f. is given by
    $$\frac{(\boldsymbol{a})_n}{(\boldsymbol{a}+\boldsymbol{b})_n} \frac{(\alpha)_n}{n!} \left(\frac{p}{1-p}\right)^n \sum_{m=0}^\infty \frac{(\alpha + n)_m}{m!} p^m (1-p)^{\alpha+n}  \pFq{K+1}{K}{-m, \boldsymbol{a} +n\mathbf{1}}{\boldsymbol{a} + \boldsymbol{b} + n\mathbf{1}}{1},$$
    where $\mathbf{1}$ is the unit vector and $\pFq{p}{q}{\mathbf{a}}{\mathbf{b}}{z}$ is the \textit{generalized hypergeometric function}.
    In particular, if $p\in (0, 1/2)$, its p.m.f. is given by
    $$\frac{(\boldsymbol{a})_n}{(\boldsymbol{a}+\boldsymbol{b})_n} \frac{(\alpha)_n}{n!} \left(\frac{p}{1-p}\right)^n \pFq{K+1}{K}{\alpha+n, \boldsymbol{a} +n\mathbf{1}}{\boldsymbol{a} + \boldsymbol{b} + n\mathbf{1}}{\frac{p}{p-1}}.$$
\end{proposition}

\begin{figure}[H]
    \centering
    \includegraphics[width=0.5\textwidth]{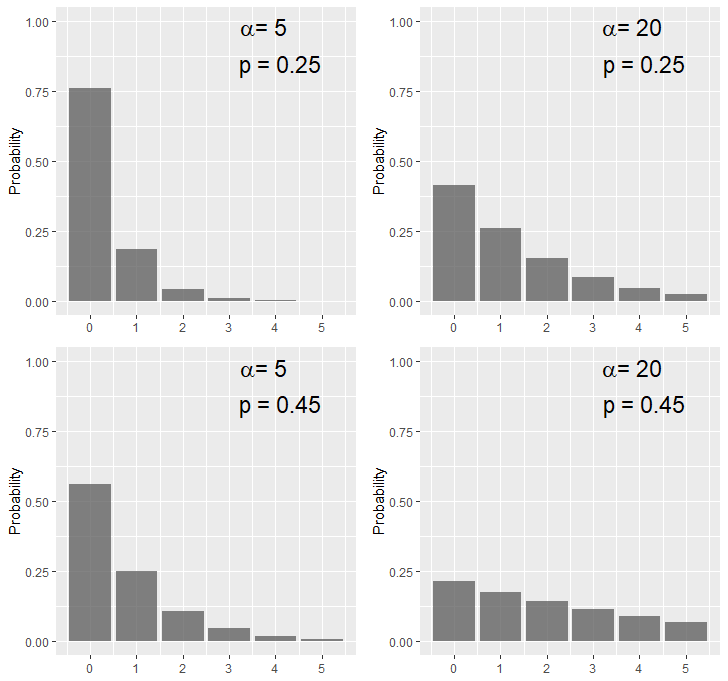}
    \caption{P.m.f. of $Y_6$ with $\theta_6 = 3$, $\theta_7 = \theta_{6,7} = 1$, $\theta_{4,5} + \theta_{8,9,10} = 2$, $\pi_{4:10} = 0.75$, $\alpha \in \{5, 20\}$ and $p = \{ 0.25, 0.45 \}$.}
    \label{fig:pmf}
\end{figure}

From Proposition \ref{prop:dispersion}, we infer that all the marginals in the example are overdispersed.
Indeed, since the negative binomial distribution is overdispersed and the tree is composed of multinomial and Dirichlet-multinomial Splittings, all the subsums and leaves are overdispersed.
For instance, the marginal distribution $Y_6$ given by \eqref{eqn:marginal} is overdispersed. 
Using Propositions \ref{prop:absorb} and \ref{prop:marginal_nb} with fixed parameters $\theta_6 = 3$, $\theta_7 = \theta_{6,7} = 1$, $\theta_{4,5} + \theta_{8,9,10} = 2$, $\pi_{4:10} = 0.75$, the p.m.f. of $Y_6$ is presented in Figure \ref{fig:pmf} for different combinations of $\alpha$ and $p$.

Now, if instead the total in the example is a fixed value $m$, i.e. a Dirac at $m$, then the nature of the Tree Pólya Splitting distribution drastically changes.
Indeed, the marginals are necessarily bounded and their dispersion may vary.
This particular case corresponds to a similar distribution as the Dirichlet-tree multinomial presented in Section \ref{subsection:Def_Tree}.
Previously, the binomial splits were absorbed by the negative binomial in accordance to Proposition \ref{prop:absorb}. 
However, because we replaced the negative binomial by a Dirac at $m$, we must find the p.m.f. of the random variable
\begin{equation}
    \label{eqn:comp_dirac}
    X \sim \bigwedge_{k = 1}^K \mathcal{P}^{[c_k]}_{n_k}\left( \theta_k, \tau_k \right) \underset{n_K}{\wedge} \mathcal{D}_m,
\end{equation}
where the parameters are given in Proposition \ref{prop:absorb}.
Using similar techniques, we prove the following.

\begin{proposition}
    \label{prop:marginal_dirac}
    Let $X$ be distributed as \eqref{eqn:comp_dirac} with $m \in \mathbb{N}$ and $M$ Pólya distributions with $c_k = 1$ and parameters $a_k, b_k \in \mathbb{R}_+$, and $K-M$ with $c_k = 0$ and parameters $\pi_k \in (0,1)$.
    Then its p.m.f. is given by
    $$\mathbb{P}(X=n) = \binom{m}{n} \frac{(\boldsymbol{a})_n}{(\boldsymbol{a}+\boldsymbol{b})_n} \gamma^n \pFq{M+1}{M}{n-m, \boldsymbol{a} +n\mathbf{1}}{\boldsymbol{a} + \boldsymbol{b} + n\mathbf{1}}{\gamma} \text{ ; } n \in \{0,\dots, m\},$$
    where $\gamma = \prod_{k=1}^{K-M} \pi_k$. 
\end{proposition}

Notice that the p.m.f. of the marginal presented in Proposition \ref{prop:marginal_dirac} generalizes those of the Dirichlet-tree multinomial and the generalized Dirichlet-multinomial. Moreover, as presented in Section 2.2.2, if the expectation of the negative binomial is equal to the Dirac distribution, i.e., $\mathrm{E}[\abs{\mathbf{Y}}] = m$, then the variance of the marginals of such a Tree Pólya Splitting, when $\mathcal{L}(\boldsymbol{\psi}) = \mathcal{NB}(\alpha, p)$, is always greater than when $\mathcal{L}(\boldsymbol{\psi}) = \mathcal{D}_m$.

\subsubsection{Factorial moments}

Factorial moments of the Pólya Splitting distribution were determined by Proposition \ref{prop:factorial_moment} to be the product of the $J$ splits.
In the same fashion, the factorial moments of the Tree Pólya splitting distribution should admit a similar product, but through its partition tree.
The following proposition shows that it is indeed the case.

\begin{proposition}
\label{prop:factorial_mom}
  For $\mathbf{Y} \sim \mathcal{TP}_{\Delta_n}(\mathfrak{T}; \boldsymbol{\theta}, \mathbf{c}) \underset{n}{\wedge} \mathcal{L}(\boldsymbol{\psi})$, $\mathbf{r} = (r_1, \dots, r_J) \in \mathbb{N}_+^J$, and by denoting $\mathbf{r}_A = (r_i)_{i\in A}$ for any $A \in \mathfrak{I} \cup \mathfrak{L}$, the factorial moments are given by
  $$\mathrm{E}[(\mathbf{Y})_{(\mathbf{r})} ] = \mu_{\abs{\mathbf{r}}} \prod_{A \in \mathfrak{I}} 
  \frac{\prod_{C \in \mathfrak{C}_A}(\theta_C)_{(\abs{\mathbf{r}_C}, c_A)}}{(\abs{\boldsymbol{\theta}_A})_{(\abs{\mathbf{r}_A}, c_A)}},$$
  where $\mu_{\abs{\mathbf{r}}}$ is the $\abs{\mathbf{r}}$-th factorial moment of $\mathcal{L}(\boldsymbol{\psi})$.
\end{proposition}

A direct corollary of Proposition \ref{prop:factorial_mom} is the univariate factorial moments of any subsum in the tree.
Here, instead of a product on all edges, the factorial moments are determined by a path from the root to the appropriate node or leaf.

\pagebreak

\begin{corollary}
\label{corollary:factorial_mom}
For $\mathbf{Y} \sim \mathcal{TP}_{\Delta_n}(\mathfrak{T}; \boldsymbol{\theta}, \mathbf{c}) \underset{n}{\wedge} \mathcal{L}(\boldsymbol{\psi})$, $r \in \mathbb{N}_+$, and any $A \in \mathfrak{I} \cup \mathfrak{L}$, the factorial moment of $\abs{\mathbf{Y}_A}$ with $A_k \in \mathrm{Path}_{A}$ is given by
  $$\mathrm{E}\left[(\abs{\mathbf{Y}_A})_{(r)}\right] = \mu_{r} \prod_{k=1}^K \frac{(\theta_{A_{k-1}})_{(r, c_{A_k})}}{(\abs{\boldsymbol{\theta}_{A_k}})_{(r, c_{A_k})}},$$
  where $\mu_{r}$ is the $r$-th factorial moment of $\mathcal{L}(\boldsymbol{\psi})$.
\end{corollary}

\noindent \textbf{Running example}

\noindent The factorial moment of $Y_6$ in the example is determined by the same path as in Figure \ref{fig:example_6}.
Corollary \ref{corollary:factorial_mom} states that for $r \in \mathbb{N}_+$,
$$\mathrm{E}\left[(Y_6)_{(r)}\right] = \mu_r \left(\pi_{4:10}^r \frac{(\theta_{6,7})_r}{(\theta_{4,5}+\theta_{6,7}+\theta_{8,9,10})_r} \frac{(\theta_{6})_r}{(\theta_{4,5}+\theta_{7})_r}\right).$$
In this case, the distribution of $\abs{\mathbf{Y}}$ is given by $\mathcal{NB}(\alpha,p)$ with $\mu_r = (\alpha)_r \left(\frac{p}{1-p} \right)^r$.
If instead the negative binomial is replaced by a Dirac at $m$, then $\mu_r = (m)_{(r)}$.
In the latter case, the factorial moment of $Y_6$ is affected only by the support of $\mathcal{D}_m$, which may be difficult to estimate in practice.
Now, let $\mathrm{E}[\abs{\mathbf{Y}}] = m$ and consider the marginals $Y_6^{\mathcal{NB}}$ and $Y_6^{\mathcal{D}}$ for the case where $\mathcal{L}(\boldsymbol{\psi})$ is given by $\mathcal{NB}(\alpha,p)$ and $\mathcal{D}_m$, respectively. 
Then, using Corollary \ref{corollary:factorial_mom} and the relation of the variance to its factorial moments, it can easily be shown that
$$\mathrm{Var}(Y_6^{\mathcal{NB}})-\mathrm{Var}(Y_6^{\mathcal{D}}) = \frac{m}{1-p}\left(\pi_{4:10}^2 \frac{(\theta_{6,7})_2}{(\theta_{4,5}+\theta_{6,7}+\theta_{8,9,10})_2} \frac{(\theta_{6})_2}{(\theta_{4,5}+\theta_{7})_2}\right) > 0.$$
Therefore, for a fixed expected value, the variance when the sum is a negative binomial is always greater than when the sum is a Dirac distribution at $m$. 
Notice that such an analysis also holds for all marginals in our example.

\subsubsection{Covariance and Correlation}

If two leaves are siblings in the partition tree, then their covariance is simply given by (\ref{eqn:covariance}) for a Pólya Splitting model with univariate distribution $\mathcal{L}(\boldsymbol{\psi})$ given by Proposition \ref{prop:marg_tree_univ}.
Therefore, the signs of covariance must be similar for any pairs of siblings. 
Hence, let us study the covariance of $Y_i$ and $Y_j$ that are not siblings. 
Using a similar argument as in Proposition \ref{prop:factorial_mom}, it can be shown that $\mathrm{Cov}(Y_i, Y_j)$ is given as follows.

\begin{proposition}
\label{prop:covariance_tree}
For $\mathbf{Y} \sim \mathcal{TP}_{\Delta_n}(\mathfrak{T}; \boldsymbol{\theta}, \mathbf{c}) \underset{n}{\wedge} \mathcal{L}(\boldsymbol{\psi})$ and marginals $Y_i, Y_j$ with $\mathscr{P}(\{i\}) \neq \mathscr{P}(\{j\})$, then there is a common ancestor node $S \in \mathfrak{I}$ with $C_i, C_j \in \mathfrak{C}_S$ such that $i \in C_i$, $j \in C_j$, and 
  \begin{equation}
      \label{eqn:cov_prop}
      \mathrm{Cov}(Y_i, Y_j) = \frac{\gamma_i \gamma_j}{\gamma_{C_i} \gamma_{C_j}} \mathrm{Cov}(\abs{\mathbf{Y}_{C_i}},\abs{\mathbf{Y}_{C_j}}),
  \end{equation}
  where $\gamma_i = \prod_{k=1}^{K} \frac{\theta_{A_{k-1}}}{\abs{\boldsymbol{\theta}_{A_k}}}$ with $A_k \in \mathrm{Path}_{i}$, and $\gamma_j$, $\gamma_{C_i}$, and $\gamma_{C_j}$ are defined analogously.
\end{proposition}

A direct consequence of this result is that the sign of covariance between $Y_i$ and $Y_j$ depends only on the ancestor node $S$. 
Therefore, as we will see in the running example, it is possible to have a covariance matrix with elements of different signs.
Precisely, at the ancestor node $S$, the sign depends only on the value of $c_S$, i.e. the type of split, and the dispersion of $\abs{\mathbf{Y}_S}$ with distribution given in Proposition \ref{prop:marg_tree_univ}.
Hence, the Tree Pólya Splitting model allows for a richer covariance structure that of the Pólya Splitting model.
Based on this result, and using Corollary \ref{corollary:factorial_mom}, we can easily develop the covariance formula in terms of all the parameters involved in the paths from the root to the leaves and the common ancestor node.

\begin{proposition}
\label{prop:covariance_tree_2}
  The covariance \eqref{eqn:cov_prop} is given by
  $$\mathrm{Cov}(Y_i, Y_j) = \gamma_i \gamma_j \left[ \left(\frac{\abs{\boldsymbol{\theta}_S}}{\abs{\boldsymbol{\theta}_S} + c_S} \right) \frac{\delta_S}{\gamma_S} \mu_2 - \mu_1^2\right],$$
  where $\delta_S = \prod_{k=1}^{K} \frac{\theta_{A_{k-1}} + c_{A_k}}{|\boldsymbol{\theta}_{A_k}| + c_{A_k}}$ with $A_k \in \mathrm{Path}_{S}$, and $\mu_{r}$ the $r$-th factorial moment of $\mathcal{L}(\boldsymbol{\psi})$.
\end{proposition}

Notice that the previous results in Proposition \ref{prop:covariance_tree} and \ref{prop:covariance_tree_2} can be adapted to any pair of subsums in the partition tree.
Another interesting result describes how the ratio of covariances is equal to a ratio of expectations.
Indeed, we directly have the following corollary.

\begin{corollary}
\label{corollary:constant}
  For $\mathbf{Y} \sim \mathcal{TP}_{\Delta_n}(\mathfrak{T}; \boldsymbol{\theta}, \mathbf{c}) \underset{n}{\wedge} \mathcal{L}(\boldsymbol{\psi})$, $A,B \in \mathfrak{I}$ such that $\mathscr{P}(A) \neq B$, $C_{A_1}, C_{A_2} \in \mathfrak{C}_A$, $C_B \in \mathfrak{C}_B$, and $B$ is not equal or descendant of $C_{A_1}$ or $C_{A_2}$,  we have
  $$\frac{\mathrm{Cov}\left(\abs{\mathbf{Y}_{C_{A_1}}},\abs{\mathbf{Y}_{C_B}}\right)}{\mathrm{Cov}\left(\abs{\mathbf{Y}_{C_{A_2}}},\abs{\mathbf{Y}_{C_B}}\right)}= \frac{\mathrm{E}\left[\abs{\mathbf{Y}_{C_{A_1}}}\right]}{\mathrm{E}\left[\abs{\mathbf{Y}_{C_{A_2}}}\right]} = \frac{\theta_{C_{A_1}}}{\theta_{C_{A_2}}}.$$
\end{corollary}

For the Pearson correlation, a similar result holds.
Indeed, since the covariances are proportional to their ancestor node, the Pearson correlation is proportional to the same node.
Using Corollary \ref{corollary:factorial_mom} for the standard deviations and Proposition \ref{prop:covariance_tree_2} for the covariance, we have the following result.

\begin{proposition}
    \label{prop:correlation_tree}
For $\mathbf{Y} \sim \mathcal{TP}_{\Delta_n}(\mathfrak{T}; \boldsymbol{\theta}, \mathbf{c}) \underset{n}{\wedge} \mathcal{L}(\boldsymbol{\psi})$, marginals $Y_i, Y_j$ with $\mathscr{P}(\{i\}) \neq \mathscr{P}(\{j\})$, then there is a ancestor node $S \in \mathfrak{I}$ such that 
  $$\mathrm{Corr}(Y_i, Y_j) = \Lambda_i \Lambda_j  \left[ \left(\frac{\abs{\boldsymbol{\theta}_S}}{\abs{\boldsymbol{\theta}_S} + c_S} \right) \frac{\delta_S}{\gamma_S} \mu_2 - \mu_1^2\right],$$
  where $$\Lambda_k = \sqrt{\frac{\gamma_k}{\delta_{k} \mu_2 + \mu_1(1-\gamma_{k} \mu_1)}} \hspace{2 pt}\text{, } k = i,j.$$ 
\end{proposition}

Observe that the correlation depends on $\Lambda_i$, $\Lambda_j$, which themselves depend on the constants $\gamma_i$, $\gamma_j$ respectively.
By Proposition \ref{prop:covariance_tree}, these constants are defined as products of values in $(0,1)$ and depend on the depth of the path from the leaf to the root, i.e. the value $K$.
Therefore, as $K$ becomes large, the constants $\gamma_i$ and $\gamma_j$ approach $0$, and so does the correlation.
In other words, if the marginals $Y_i$ and $Y_j$ are distant in the tree, their correlation tends to $0$.\\

\noindent \textbf{Running example}

\begin{figure}[ht]
    \centering
    \resizebox{0.5\textwidth}{0.2\textheight}{
\begin{tikzpicture}[scale = 0.9]
    \node{$\mathcal{NB}(\alpha,p)$}
    child[line width= 3 pt, dashed]{ node {$\mathcal{M}$}[sibling distance = 3.5cm]
            child[thick, solid] {node[label={[xshift = 1.4cm,yshift= 0.4cm] $\pi_{1,2}$},] {$\mathcal{DM}$}[sibling distance = 2cm]
                child{node[label={[xshift = 0.1cm,yshift= 0.2cm] $\theta_{1}$}]{$Y_1$}}
                child{node[label={[xshift = -0.1cm,yshift= 0.2cm] $\theta_{2}$}]{$Y_2$}}
            }
            child[thick, solid] {node[label={[xshift = -0.4cm,yshift= 0.1cm] $\pi_{3}$}] {$Y_3$}
            }
            child{node[label={[xshift = -1.4cm,yshift= 0.4cm] $\pi_{4:10}$}, label={[xshift = 0.8cm,yshift= 0cm] $S$}, draw, solid, line width = 2pt] {$\mathcal{DM}$}[sibling distance = 3cm]
      child[thick, solid] {node[label={[xshift = 0.8cm,yshift= 0.2cm] $\theta_{4,5}$}] {$\mathcal{M}$}[sibling distance = 1.5cm]
                child{node[label={[xshift = -0.1cm,yshift= 0.2cm] $\pi_{4}$}]{$Y_4$}}
                child[thick, solid]{node[label={[xshift = 0.1cm,yshift= 0.2cm] $\pi_{5}$}]{$Y_5$}}
      }
      child[line width = 3pt, solid] {node[label={[xshift = -0.45cm,yshift= 0cm] $\theta_{6,7}$}] {$\mathcal{DM}$}[sibling distance = 1.5cm]
                child{node[label={[xshift = -0.1cm,yshift= 0.2cm] $\theta_{6}$}]{$Y_6$}}
                child[thick]{node[label={[xshift = 0.1cm,yshift= 0.2cm] $\theta_{7}$}]{$Y_7$}}
      }
      child[line width = 3pt, solid] {node[label={[xshift = -0.7cm,yshift= 0.2cm] $\theta_{8,9,10}$}] {$\mathcal{DM}$}[sibling distance = 1.5cm]
                child[thick]{node[label={[xshift = -0.1cm,yshift= 0.2cm] $\theta_{8}$}]{$Y_8$}}
                child{node[label={[xshift = 0.1cm,yshift= 0.2cm] $\theta_{9,10}$}]{$\mathcal{M}$}
                    child{node[label={[xshift = 0cm,yshift= 0.2cm] $\pi_{9}$}]{$Y_9$}}
                    child[thick]{node[label={[xshift = 0cm,yshift= 0.2cm] $\pi_{10}$}]{$Y_{10}$}}            
                }
      }
    }};
\end{tikzpicture}
}
    \caption{Covariance between $Y_6$ and $Y_9$ in the running example.}
    \label{fig:covariance}
\end{figure}
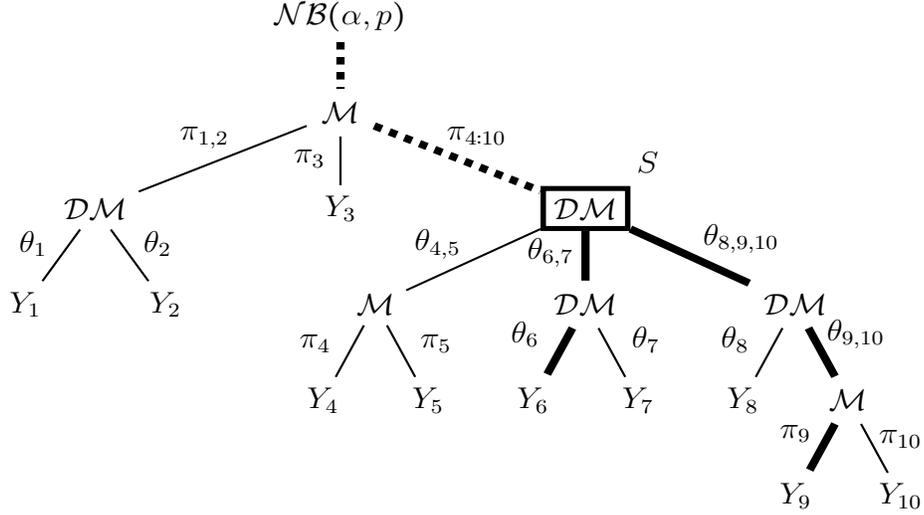

\noindent Let us study the covariance for the running example.
Since $\mathcal{L}(\boldsymbol{\psi}) = \mathcal{NB}(\alpha,p)$, the overdispersion is preserved throughout the tree by Proposition \ref{prop:dispersion}. 
Therefore, it should be possible to have different covariance signs.
We begin by calculating the covariance between $Y_6$ and $Y_9$.
As presented in Figure \ref{fig:covariance}, the ancestor node is given by $S = \{4, \dots, 10\}$ with $C_6 = \{6,7\}$ and $C_9 = \{8,9,10\}$.
By Propositions \ref{prop:covariance_tree} and \ref{prop:covariance_tree_2}, we have for $\boldsymbol{\theta}_S = \{\theta_{4,5}, \theta_{6,7}, \theta_{8,9,10}\}$
$$\gamma_S = \delta_S =  \pi_{4:10}, \hspace{20 pt} \gamma_6 = \pi_{4:10} \cdot \frac{\theta_{6,7}}{\abs{\boldsymbol{\theta}_S}} \cdot \frac{\theta_{6}}{\theta_{6} + \theta_{7}}, \hspace{20 pt} \gamma_9 = \pi_{4:10} \cdot \frac{\theta_{8,9,10}} {\abs{\boldsymbol{\theta}_S}}\cdot \frac{\theta_{9,10}}{\theta_{8} + \theta_{9,10}} \cdot \pi_9.$$
Using the identity $\mu_r = (\alpha)_r p^r/(1-p)^r$ for a negative binomial distribution, the covariance is given by
\begin{align}
\label{eqn:cov_example}
    \mathrm{Cov}(Y_6, Y_9) &= \left[\frac{\theta_{6}}{\theta_{6} + \theta_{7}}\right] \left[\pi_9 \frac{\theta_{9,10}}{\theta_{8}+\theta_{9,10}}\right] \mathrm{Cov}(Y_6+Y_7, Y_8+Y_9+Y_{10}) \nonumber\\
    &= \alpha \left(\frac{p}{1-p} \right)^2 \left[(\alpha+1)\left(\frac{\abs{\boldsymbol{\theta}_S}}{\abs{\boldsymbol{\theta}_S} + 1}\right) -  \alpha \right] \gamma_6  \gamma_9 \nonumber\\
    &= \alpha \left(\frac{p}{1-p} \right)^2 \left(\frac{\abs{\boldsymbol{\theta}_S} - \alpha}{\abs{\boldsymbol{\theta}_S} + 1}\right) \gamma_6  \gamma_9. 
\end{align}
Equation \eqref{eqn:cov_example} indicates that the covariance can either be negative, positive or null whether $\alpha$ is greater, smaller or equal to $\abs{\boldsymbol{\theta}_S}$ respectively.
In fact, using Theorem 6 of \cite{PEYHARDI2021104677}, the distribution of $\abs{\mathbf{Y}_S}$ is given by $$\abs{\mathbf{Y}_S} \sim \mathcal{NB}\left(\alpha, \frac{p \pi_{4:10}}{1-p+p\pi_{4:10}}\right),$$
and since the negative binomial is the only distribution that can induce independence for Dirichlet-multinomial Splitting, $Y_6$ and $Y_9$ are independent if and only if $\alpha = \abs{\boldsymbol{\theta}_S}$.
Therefore, a null correlation truly indicates independence in this case.
All of these dependence relationships also hold true for any pair of random variables between the subsets $\{Y_4, Y_5\}$, $\{Y_6, Y_7\}$ and $\{Y_8, Y_9, Y_{10}\}$ since their common ancestor node is still $S$.

For the sake of further investigation, suppose that $\alpha = \abs{\boldsymbol{\theta}_S}$, but $\alpha > \theta_1 + \theta_2$.
Given that the distribution of $Y_1 + Y_2$ is also negative binomial, a similar argument leads us to conclude that $Y_1$ and $Y_2$ are negatively correlated, whereas $Y_6$ and $Y_9$ are independent.
For a concrete example, let us fix the model parameters as:

\begin{minipage}{.3\linewidth}
\small
\begin{align*}
    \alpha &= 10\\
    p&=0.95\\
    \pi_{1,2} = \pi_9 &= 0.3\\
    \pi_3 &= 0.1\\
    \theta_1 = \theta_2 &= 1.5
\end{align*}
\end{minipage}
\hfill
\begin{minipage}{.3\linewidth}
\begin{align*}
  \pi_{4:10} &= 0.6\\
  \theta_{4,5} &= 3\\
  \theta_{6,7} = \theta_{8,9.10} &= 3.5\\
  \pi_4 = \pi_5 &= 0.5
\end{align*}
\end{minipage}
\hfill
\begin{minipage}{.3\linewidth}
\begin{align}
  \label{parameters}
  \theta_6 &= 0.8 \nonumber\\
  \theta_7 = \theta_8 &= 1\\
  \theta_{9,10} &= 2.5\nonumber\\
  \pi_{10} &= 0.7 \hspace{4pt}.\nonumber
\end{align}
\end{minipage}

\normalsize

\begin{figure*}[ht]
    \centering
    \begin{subfigure}[b]{0.5\textwidth}
        \centering
        \includegraphics[width=0.6\textwidth]{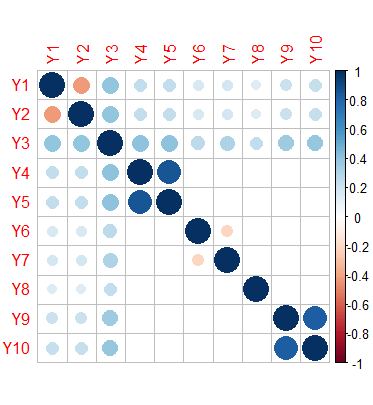}
        \caption{$\mathcal{L}(\boldsymbol{\psi}) = \mathcal{NB}(\alpha, p)$}
        \label{fig:corr_example}
    \end{subfigure}%
    ~ 
    \begin{subfigure}[b]{0.5\textwidth}
        \centering
        \includegraphics[width=0.6\textwidth]{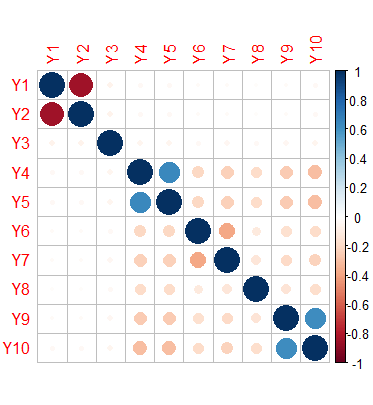}
        \caption{$\mathcal{L}(\boldsymbol{\psi}) = \mathcal{D}_{100}$}
        \label{fig:corr_example2}
    \end{subfigure}
    \caption{Correlation plot of the running example with parameters \eqref{parameters} and different $\mathcal{L}(\boldsymbol{\psi})$.}
\end{figure*}

\noindent The correlation can be computed using Proposition \ref{prop:correlation_tree} and be summarized in Figure \ref{fig:corr_example}.
To conclude this section, let us study the same example, but with $\mathcal{L}(\boldsymbol{\psi})$ as a Dirac at $m$. 
With this simple modification, the covariance at the root must be negative since $\mathcal{D}_m$ is underdispersed.
Furthermore, it is possible that all marginals remain underdispersed by Proposition \ref{prop:dispersion}. 
If the hypergeometric distribution is introduced in this example, then it would be possible to have different signs of covariance.
For sake of comparison, suppose $m = 100$ and the tree parameters take the same values as before. 
Again, using Proposition \ref{prop:correlation_tree}, the correlations can be calculated and they are presented in Figure \ref{fig:corr_example2}.\vspace{-3 pt}

\pagebreak

\subsubsection{Log-likelihood decomposition}
\label{subsubsection:loglik}

As obtained by \cite{peyhardi2016partitioned}, a direct consequence of Definition \ref{def:Tree_Pólya} is the decomposition of the log-likelihood with respect to the partition tree.
Indeed, the log-likelihood of \eqref{eqn:tree_polya} is given by

\begin{equation}
    \label{eqn:loglikelihood}
    \log\big[\mathbb{P}(\abs{\mathbf{Y}} = n)\big] + \sum_{A \in \mathfrak{I}} \left(  \sum_{C \in \mathfrak{C}_A} \log (\theta_C)_{(n_C, c_A)} - \log (\abs{\boldsymbol{\theta}_A})_{(n_A,c_A)}  \right)  + \mathrm{constant.}
\end{equation}
Therefore, the maximum likelihood estimators (MLE), if they exist, of the whole Tree Pólya Splitting can be obtained by combining the MLE of $\boldsymbol{\psi}$ based on $\abs{\mathbf{Y}} \sim \mathcal{L}(\boldsymbol{\psi})$ and each Pólya on $\mathfrak{T}$ separately.
This divide-and-conquer approach greatly simplifies inference and potentially model selection as well.
Indeed, the AIC of the whole model is simply the sum of those at each node.

\section{Analysis of a Trichoptera dataset}
\label{section:application}

In this section, we fit a Tree Pólya Splitting distribution with a fixed partition tree $\mathfrak{T}$ to the Trichoptera dataset provided by \cite{usseglio1987influence} and compare it to other splittings models with a multinomial, Dirichlet-multinomial and generalized Dirichlet-multinomial, as well as to the multivariate Poisson-lognormal distribution \citep{AITCHISON, PLNmodels}.
We also investigate how the data can be used to find the partition tree $\mathfrak{T}$.
Specifically, we propose an algorithm that searches for an appropriate $\mathfrak{T}$ to fit the Tree Pólya Splitting using the identified structure.
The Trichoptera dataset consists of $J=17$ species' abundances and  $9$ covariates collected from $n=49$ sites between 1959 and 1960.
For the sake of simplicity, no covariates are used in this application.
However, it is important to note that incorporating them into the model is feasible.
As presented in \citep{PEYHARDI2021104677}, the log-likelihood decomposition of the Pólya Splitting remains valid, which allows the model to include both a univariate regression component and a Pólya regression component.
The latter can be implemented using methods such as those presented in \cite{Zhang2017}.
For the Pólya Tree Splitting, a similar approach can be made.

For the total distribution, the data exhibits overdispersion.
Likewise, all species except three exhibit empirical overdispersion. 
Specifically, the species \textit{Che}, \textit{Hyc}, and \textit{Hys} appear to be underdispersed.
As previously explained, the log-likelihood of the Tree Pólya can be decomposed with respect to the partition tree.
Therefore, the parameters can be estimated step-by-step starting by the distribution of the total $\mathcal{L}(\boldsymbol{\psi})$, and then each Pólya Splitting in the tree.
We fix $\mathcal{L}(\boldsymbol{\psi})$ to be a negative binomial since the data are overdispersed, and we estimate its parameters by MLE using the \texttt{R} package \texttt{MASS} \citep{MASS}.
This yields a negative binomial with parameters $\alpha = 0.478$, $p = 0.997$ and an AIC of $575.016$.
Since this distribution has an unbounded support, either a multinomial or Dirichlet-multinomial can be adjusted at each internal node of $\mathfrak{T}$ using the \texttt{R} package \texttt{MGLM} \citep{Zhang2017}.

Now, the partition tree $\mathfrak{T}$ remains to be fixed.
A natural approach is to use evolutionary information concerning the Trichoptera.
Indeed, since the data consists of $J=17$ different species, it is possible to regroup them into respective families.
Precisely, we can divide the $17$ species into $5$ different groups based on the information provided by \cite{usseglio1987influence}. 
With this structure, we adjust either a multinomial or a Dirichlet-multinomial at each node based on the AIC.
However, should the Dirichlet-multinomial parameter estimates fail to converge, we opt for a multinomial distribution.
Indeed, when $\mathcal{DM}(\boldsymbol{\theta})$ is such that each $\theta_j \to \infty$ and $\theta_j/\abs{\boldsymbol{\theta}} \to \pi_j \in (0,1)$ for all $j \in \{1,\dots,J\}$, then the Dirichlet-multinomial converges to a multinomial with parameters $\boldsymbol{\pi} = (\pi_1, \dots, \pi_J)$.\vspace{-5 pt}

\begin{figure}[H]
    \centering
    \resizebox{0.45\textwidth}{0.16\textheight}{
\begin{tikzpicture}[scale = 0.7]
    \node{$\mathcal{NB}$}
        child{ node{$\mathcal{DM}$}[sibling distance = 10cm]
            child{node{$\mathcal{DM}$}[sibling distance = 4.8cm]
                    child{node{$\mathcal{M}$}[sibling distance = 1.2cm]
                        child{node{Che}}
                        child{node{Hyc}}
                        child{node{Hym}}
                        child{node{Hys}}                    
                    }
                    child{node{$\mathcal{DM}$}[sibling distance = 1.2cm]
                        child{node{Ath}}
                        child{node{Cea}}
                        child{node{Ced}}
                        child{node{Set}}                     
                    }
                    child{node{$\mathcal{DM}$}[sibling distance = 1.2cm]
                        child{node{Aga}}
                        child{node{Glo}}}
                    child{node{$\mathcal{DM}$}[sibling distance = 1.2cm]
                        child{node{All}}
                        child{node{Han}}
                        child{node{Hfo}}
                        child{node{Hsp}}
                        child{node{Hve}}
                        child{node{Sta}}
                    }
            }
            child{node{Psy}}
        };
\end{tikzpicture}
}
\caption{Tree Pólya Splitting fitted to the Trichoptera dataset with a fixed partition tree.}
\label{fig:tri_1}
\end{figure}
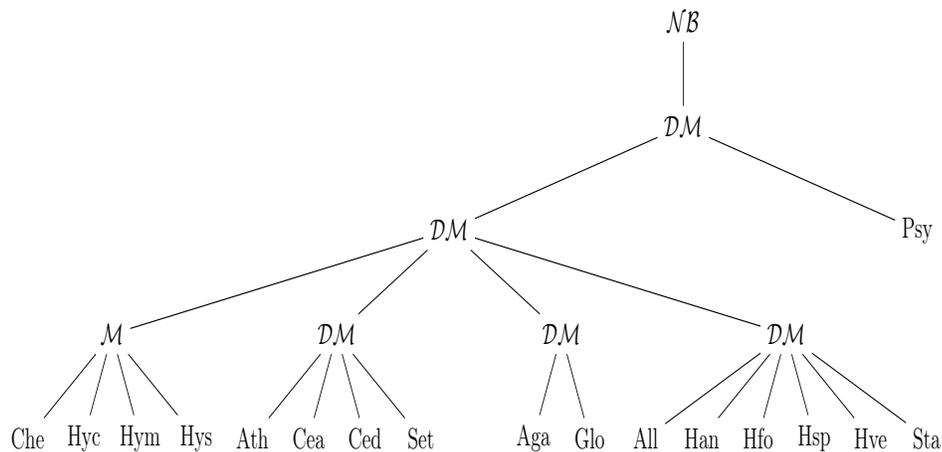

With these criteria, we have the Tree Pólya model presented in Figure \ref{fig:tri_1}.
By decomposition of the AIC with respect to the partition tree, the adjusted Tree Pólya Splitting has a AIC of $2465.85$.
If we adjust the data to a multivariate Poisson-lognormal using the \texttt{R} package \texttt{PLNmodels} \citep{PLNmodels}, we obtain an AIC of $2599.63$.
For this application, $170$ parameters are needed for the Poisson-lognormal compared to $23$ parameters for Tree Pólya Splitting.
Therefore, the proposed model is simpler and equally adequate according to the AIC score.
If, instead of the partition tree in Figure \ref{fig:tri_1}, we use the structure of a multinomial, Dirichlet-multinomial or generalized Dirichlet-multinomial, then we would obtain AICs of $6362.20$, $2494.87$, and $2460.70$ respectively.
Notice the latter is slightly better than the proposed tree structure.
In order to find a better partition tree, we must build it using the data.
Since an exhaustive search of every possible tree is a formidable task, a searching algorithm must be used to efficiently find a suitable structure.
Again, since the log-likelihood of the whole model is the sum of log-likelihoods at each internal node, the tree can be simply built step-by-step by summing the AIC.
In the following, we present an alternative approach to constructing $\mathfrak{T}$.

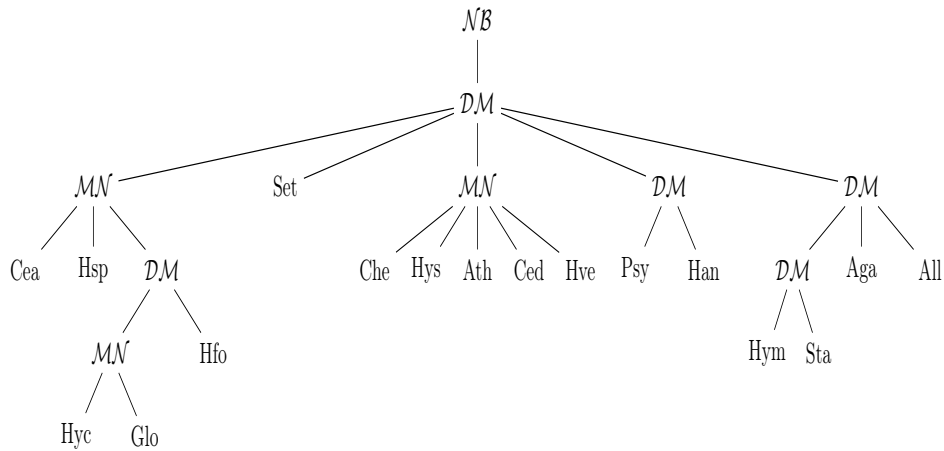
\begin{figure}[ht]
    \centering
    \resizebox{0.5\textwidth}{0.18\textheight}{
\begin{tikzpicture}[scale = 0.6]
    \node{$\mathcal{NB}$}
    child{ node{$\mathcal{DM}$}[sibling distance = 5.6cm]
        child{node{$\mathcal{M}$}[sibling distance = 2cm]
                child{node{Cea}}
                child{node{Hsp}}
                child{node{$\mathcal{DM}$}[sibling distance = 3cm]
                    child{node{$\mathcal{M}$}[sibling distance = 2cm]
                         child{node{Hyc}}
                         child{node{Glo}}                        
                    }
                    child{node{Hfo}}
                }
        }
        child{node{Set}}
        child{node{$\mathcal{M}$}[sibling distance = 1.5cm]
                child{node{Che}}
                child{node{Hys}}
                child{node{Ath}}
                child{node{Ced}}
                child{node{Hve}}
        }        
        child{node{$\mathcal{DM}$}[sibling distance = 2cm]
                child{node{Psy}}
                child{node{Han}}
        }
        child{node{$\mathcal{DM}$}[sibling distance = 2cm]
                child{node{$\mathcal{DM}$}[sibling distance = 1.5cm]
                    child{node{Hym}}
                    child{node{Sta}}
                }
                child{node{Aga}}
                child{node{All}}
        }
    };
\end{tikzpicture}
}
\caption{Tree Pólya Splitting fitted to the Trichoptera dataset with a partition tree search.}
\label{fig:tri_2}
\end{figure}

The algorithm begins with $\mathfrak{T}$ as a rooted tree, where the only internal node is $\Omega$, and $\mathfrak{C}_\Omega = \mathfrak{L}$.
For such a tree, a Dirichlet-multinomial is fitted.
As the first step, the algorithm tests whether adding an internal node improves the model.
Precisely, for some leaves $i$ and $j$, it evaluates whether the AIC improves if a new internal node $A = \{i,j\}$ is added such that $\mathfrak{C}_A = \{\{i\},\{j\}\}$, $\mathfrak{C}_\Omega = \{A, \mathfrak{L}_{-\{i,j\}}\}$, and a Dirichlet-multinomial is fitted at each internal node.
The algorithm examines every possible combinations of $\{i,j\}$, and if none improves the model, it stops and returns the initial Dirichlet-multinomial model.
Otherwise, the combination that yields the best AIC improvement is selected.

After selecting an internal node $A = \{i,j\}$, the algorithm tests whether there exists a leaf $k \in \mathfrak{L}_{-\{i,j\}}$, connected to $\Omega$, such that the model improve if this leaf is transferred to $A$.
In other words, it tests whether the model improves when a Tree Pólya Splitting with Dirichlet-multinomial splits is fitted where $A = \{i,j,k\}$, $\mathfrak{C}_A = \{\{i\},\{j\},\{k\}\}$ and $\mathfrak{C}_\Omega = \{A, \mathfrak{L}_{-\{i,j,k\}}\}$.
The algorithm tests for all possible leaves $k$ and keeps only the best in term of AIC.
This leaf transfer step is repeated as long as possible and continues to improve the model.

\begin{figure}[ht]
    \centering
    \resizebox{0.6\textwidth}{0.15\textheight}{
\begin{tabular}{c|c|c}
    \begin{tikzpicture}[scale = 0.9]
    \node{$\mathcal{DM}$}[sibling distance = 0.8cm]
            child{node{$Y_1$}}
            child{node{$Y_2$}}
            child{node{$Y_3$}}
            child{node{$Y_4$}}
            child{node{$Y_5$}}
            child{node{$Y_6$}}
            child{node{$Y_7$}};
            \node at (0,-2.8) {Iteration 0};
\end{tikzpicture}
&
    \begin{tikzpicture}[scale = 0.9]
    \node{$\mathcal{DM}$}[sibling distance = 0.8cm]
            child{node{$\mathcal{DM}$}[sibling distance = 0.8cm]
                        child{node{$Y_1$}}
                        child{node{$Y_3$}}     
            }
            child{node{$Y_2$}}
            child{node{$Y_4$}}
            child{node{$Y_5$}}
            child{node{$Y_6$}}
            child{node{$Y_7$}};
            \node at (0,-3.8) {Iteration 1};
\end{tikzpicture}
&
\begin{tikzpicture}[scale = 0.9]
    \node{$\mathcal{DM}$}[sibling distance = 0.8cm]
            child{node{$\mathcal{DM}$}[sibling distance = 0.8cm]
                        child{node{$Y_1$}}
                        child{node{$Y_3$}}
                        child{node{$Y_4$}}
            }
            child{node{$Y_2$}}
            child{node{$Y_5$}}
            child{node{$Y_6$}}
            child{node{$Y_7$}};
            \node at (0,-3.8) {Iteration 2};
\end{tikzpicture}\\
\hline
\begin{tikzpicture}[scale = 0.9]
    \node{$\mathcal{DM}$}[sibling distance = 0.9cm]
            child{node{$\mathcal{DM}$}[sibling distance = 0.8cm]
                        child{node{$Y_1$}}
                        child{node{$Y_3$}}
                        child{node{$Y_4$}}
            }
            child{node{$Y_2$}}
            child{node{$Y_6$}}
            child{node{$\mathcal{DM}$}[sibling distance = 0.8cm]
                        child{node{$Y_5$}}
                        child{node{$Y_7$}}
            };
            \node at (0,-3.8) {Iteration 3};
\end{tikzpicture}
&
\begin{tikzpicture}[scale = 0.9]
    \node{$\mathcal{DM}$}[sibling distance = 1cm]
            child{node{$\mathcal{DM}$}[sibling distance = 0.6cm]
                        child{node{$Y_1$}}
                        child{node{$Y_3$}}
                        child{node{$Y_4$}}
            }
            child{node{$Y_6$}}
            child{node{$\mathcal{DM}$}[sibling distance = 0.6cm]
                        child{node{$Y_2$}}
                        child{node{$Y_5$}}
                        child{node{$Y_7$}}
            };
            \node at (0,-4.3) {Iteration 5};
\end{tikzpicture}
&
\begin{tikzpicture}[scale = 0.9]
    \node{$\mathcal{DM}$}[sibling distance = 1cm]
            child{node{$\mathcal{DM}$}[sibling distance = 0.8cm]
                        child{node{$Y_1$}}
                        child{node{$\mathcal{DM}$}[sibling distance = 0.6cm]
                            child{node{$Y_3$}}
                            child{node{$Y_4$}}                        
                        }
            }
            child{node{$Y_6$}}
            child{node{$\mathcal{DM}$}[sibling distance = 0.6cm]
                        child{node{$Y_2$}}
                        child{node{$Y_5$}}
                        child{node{$Y_7$}}
            };
            \node at (0,-5) {Iteration 6};
\end{tikzpicture}
\end{tabular}
}
\caption{Example of the search algorithm for $\mathbf{Y} = (Y_1, \dots, Y_7)$, where each iteration corresponds to the steps of internal node creation or leaf transfer.}
\label{fig:algo}
\end{figure}

Once all possible transfers to $A$ are made, the tree $\mathfrak{T}$ is updated such that $\mathfrak{C}_\Omega =\{A, \mathfrak{L}_{-A}\}$ where $\mathfrak{L}_{-A}$ are the remaining leaves not used in the construction of $A$.
Then, the algorithm test whether all previous steps can be repeated with another node $B$, using the leaves $\mathfrak{L}_{-A}$.
This set of instruction is repeated until none remain available or if the AIC measure does not improve.
Finally, the entire procedure is repeated for all newly created internal nodes (i.e., $A, B, \dots$), and the search ultimately stops when no further improvements are achievable.
See Figure \ref{fig:algo} for an example of this search algorithm for a discrete vector $\mathbf{Y} = (Y_1, \dots, Y_7)$.

While this search algorithm provides a suitable tree structure, it may not yield the optimal one.
Firstly, it is possible that a tree structure exists with a better AIC score, which cannot be achieved by this algorithm. 
Secondly, each node in the tree has a Dirichlet-multinomial distribution.
The model may be improved by changing some nodes to a multinomial.
For each node, we test whether the AIC improves with a multinomial, or if the parameters of some Dirichlet-multinomial diverges. 
Thanks to this approach, the resulting model is presented in Figure \ref{fig:tri_2}.

\begin{table}[ht]
\centering
\begin{adjustbox}{width=0.5\textwidth}
\small
\begin{tabular}{|c|c|c|}
\hline
 Model & Nb. Parameters &  AIC\\ 
\hline
    Partition tree search (Figure \ref{fig:tri_2}) & 23 & 2380.77\\
\hline
    Generalized Dirichlet-multinomial Splitting & 34 & 2460.70\\
\hline
    Fixed partition tree (Figure \ref{fig:tri_1}) & 23 & 2465.85\\
\hline
    Dirichlet-multinomial Splitting & 19 & 2494.87\\
\hline
    Multivariate Poisson-lognormal & 170 & 2599.63\\
\hline
    Multinomial Splitting & 18 & 6362.20\\
\hline
\end{tabular}
\end{adjustbox}
\caption{Fitted models to the Trichoptera dataset.}
\label{Tab:tri}
\end{table}

\pagebreak

\begin{figure}[ht]
    \centering
    \includegraphics[width=0.45\textwidth]{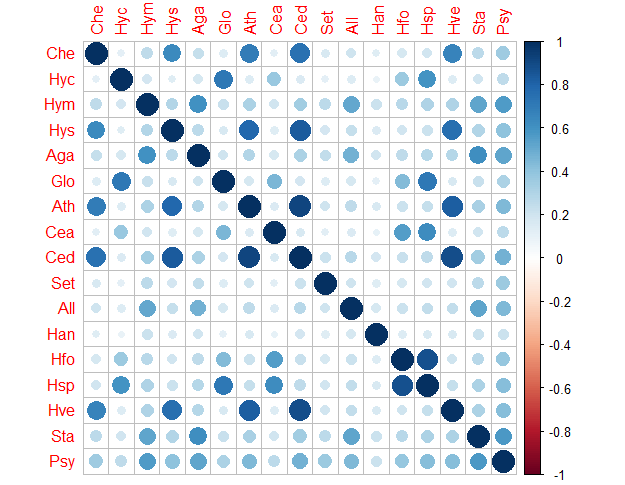}
    \caption{Correlation plot of the Tree Pólya Splitting fitted to the Trichoptera dataset with a partition tree search.}
    \label{fig:corr_tri}
\end{figure}

Using the same negative binomial distribution for $\mathcal{L}(\boldsymbol{\psi})$, the AIC of this latest model is $2380.77$, and still only needs $23$ parameters.
Each fitted model is presented in Table \ref{Tab:tri}.
Finally, the estimated correlations can easily be calculated by combining Proposition \ref{prop:correlation_tree} and the estimated parameters.
See Figure \ref{fig:corr_tri}.

\section*{Discussion and perspectives}

The simplicity and versatility of the Pólya Splitting model proposed by \cite{JONES201983} and \cite{PEYHARDI2021104677} has been expanded in this work.
The Tree Pólya Splitting model provides not only a generalization of the latter, but also allows for more diverse correlation structure. 
The partition tree of the Tree Pólya Splitting provides a convenient model for inference and a simpler parametrization. 
This paper provides the basis of Tree Pólya Splitting models and further issues remain to be explored.
Initially, \cite{PEYHARDI201759} studied the probabilistic graphical model associated to each type of Pólya Splitting.
They proved that their graphs are either complete, meaning there is an edge between all nodes, or empty, meaning no edges are present.
Specifically, \cite{Peyhardi2022} showed that the probabilistic graphical model of a Pólya Splitting distribution is empty if and only if $\mathcal{L}(\boldsymbol{\psi})=\mathcal{PS}^{[c]}(\theta, \alpha)$, while being complete otherwise, and they extended this result for a broader class of Splitting distributions, which utilizes the quasi-Pólya distribution \citep[e.g.][]{Janardan1977AGO}.
Given that the Tree Pólya Splitting model exhibits various dependence structures, it should lead to more complex graphs than those that are complete or empty, thus offering interesting avenues to the problem of learning graphical models with discrete variables.

Zero-inflation for multivariate count data would be an interesting avenue to explore as well.
Several models, including those presented by \cite{Multi_Zero_Poisson}, \cite{Multi_Zero_Poisson_2}, and \cite{Zero_logistic}, have been proposed.
Similarly, \cite{Tang&Chen2018} provide a solution to zero-inflation for the generalized Dirichlet-multinomial model.
In their work, each combination of zero-inflation is made possible thanks to the underlying binary cascade tree.
Precisely, they use the zero-inflated beta-binomial distribution at each internal node, which can be interpreted as zero-inflation at all the left leaves of the tree.
A generalization of this idea could be made for Tree Pólya Splitting, but zero-inflation could instead be modelled on any branch of the tree. 
This particular idea is explore by \cite{Fabrice} for the binary Dirichlet-tree multinomial model.

Finally, the analysis of extreme values for the Pólya and Tree Pólya Splittings could lead to interesting results.  
Indeed, the field of extreme value theory provides a wide range of results for multivariate continuous distributions, but is lacking in terms of multivariate discrete distributions.
\cite{feidt2010asymptotics} attempted to provide some answers to this problem using extreme copulas. 
\cite{valiquette} also explored this question, but for univariate Poisson mixtures. 
Given that the Tree Pólya Splitting model combines a univariate discrete distribution with a tree singular distribution, integrating both their results in this model could offer valuable insights into the challenge of modelling multivariate discrete extreme values.

\section*{Acknowledgments}

This research was supported by the GAMBAS project funded by the French National Research Agency (ANR-18-CE02-0025) and the French national programme LEFE/INSU.
Éric Marchand's research is supported in part by the Natural Sciences and Engineering Research Council of Canada.
We extend our gratitude to Stéphane Girard, Professor Stéphane Robin, and Professor Klaus Herrmann for their insightful comments on the manuscript.

\section*{Appendix - Proofs of Propositions}

\subsection*{Proof of Proposition \ref{prop:factorial_moment}}

As presented in equation 40.19 of \cite{johnson1997discrete}, the Chu-Vandermonde identity gives us the following:
$$\mathrm{E}[(\mathbf{Y})_{(\mathbf{r})}  \mid \hspace{2pt} \abs{\mathbf{Y}} = n] = \frac{(n)_{(\abs{\mathbf{r}})}}{(\abs{\boldsymbol{\theta}})_{(\abs{\mathbf{r}}, c)}} \prod_{j=1}^J (\theta_j)_{(r_j,c)}.$$
We conclude by iterated expectations.\qed

\subsection*{Proof of Proposition \ref{prop:dispersion}}

We analyze the ratio 
$$R := \frac{\mathrm{Var}[Y_j]}{\mathrm{E}[Y_j]} -1 = \frac{\mu_{2}(\theta_j + c)}{\mu_{1}(\abs{\boldsymbol{\theta}} + c)} - \frac{\mu_1 \theta_j}{\abs{\boldsymbol{\theta}}},$$
and only need to provide an inequality on the right-hand side.
For $c = 0$, $R$ is equal to
$\frac{(\mu_{2} - \mu_1^2) \theta_j}{\mu_{1} \abs{\boldsymbol{\theta}}},$
which is zero, positive or negative whether $\mathcal{L}(\boldsymbol{\psi})$ has null, over, or under dispersion, respectively. 
For $c = 1$, suppose $\mu_2 - \mu_1^2 \geq 0$.
Then $R$ is such that 
$$R \propto \mu_{2}(\theta_j + 1)\abs{\boldsymbol{\theta}}- \mu_{1}^2\theta_j (\abs{\boldsymbol{\theta}}+1) = (\mu_2 - \mu_1^2)\theta_j(\abs{\boldsymbol{\theta}}+1) + \mu_2\abs{\boldsymbol{\theta}_{-j}},$$
which is positive. A similar argument for $c=-1$ shows that $R < 0$ when $\mu_2 - \mu_1^2 \leq 0$.\qed

\subsection*{Proof of Proposition \ref{prop:correlation}}

By Proposition \ref{prop:factorial_moment}, we have for any $\ell \in \{i,j\}$ and $r \in \mathbb{N}$ that $\mu_r = \frac{(\abs{\boldsymbol{\theta}})_{(r,c)}}{(\theta_\ell)_{(r,c)}} \mathrm{E}[(Y_\ell)_{(r)}]$.
Using this identity in \eqref{eqn:covariance}, we have 
\begin{align*}
    \mathrm{Cov}(Y_i, Y_j) &= \frac{\theta_j}{\theta_i + c} \mathrm{Var}[Y_i](1-M_i)\\
    &= \frac{\theta_i}{\theta_j + c} \mathrm{Var}[Y_j](1-M_j),
\end{align*}
implying that $\mathrm{sgn}(1-M_i) = \mathrm{sgn}(1-M_j)$.
Suppose without lost of generality that $M_i \neq 1$, then the variance of $Y_i$ is such that $$\mathrm{Var}[Y_i] = \frac{(\theta_i + c)\mathrm{Cov}(Y_i, Y_j)}{\theta_j(1-M_i)}.$$
A similar equality for $\mathrm{Var}[Y_j]$ can be obtained. 
By definition of the Pearson correlation coefficient, we then have \eqref{eqn:corr}.\qed

\subsection*{Proof of Proposition \ref{prop:bound_corr}}

It is sufficient to show this result when the parameters of $\mathbf{Y}\sim \mathcal{DM}_{\Delta_n}(\boldsymbol{\theta}) \underset{n}{\wedge} \mathcal{L}(\boldsymbol{\psi})$ allow the covariance to be positive. 
Under this hypothesis, $1-M_i$ and $1-M_j$ must be positive.
Therefore $(1-M_i)(1-M_j) < 1$, and we can conclude by applying this inequality to \eqref{eqn:corr}.\qed

\subsection*{Proof of Propositions \ref{prop:marg_tree_univ} and \ref{prop:marg_tree_multi}}

By definition, a Tree Pólya model has Pólya Splitting at each internal node. 
In particular, for the root $\Omega$, each marginal is a subsum $\abs{\mathbf{Y}_{C_j}}$, with $C_j \in \mathfrak{C}(\Omega)$, such that their distribution is given by \eqref{eqn:marg} with $$\abs{\mathbf{Y}_{C_j}} \sim \mathcal{P}_{\Delta_n}^{[c_{\Omega}]}\left(\theta_{C_j}, \abs{\boldsymbol{\theta}_{\Omega \setminus C_j}} \right) \underset{n}{\wedge} \mathcal{L}(\boldsymbol{\psi}).$$
For each $j \in \{1, \dots, J_\Omega\}$, we have an induced univariate distribution that defines new roots in the tree. 
Iterating this process, we conclude that the three propositions follow.\qed

\subsection*{Proof of Proposition \ref{prop:absorb}}

From Theorem $6$ of \cite{PEYHARDI2021104677}, $\mathcal{B}_n(\pi) \underset{n}{\wedge} \mathcal{NB}(\alpha,p) = \mathcal{NB}\left(\alpha, \frac{p \pi}{1- p(1-\pi)}\right)$. 
By iterating this composition, we can easily show that 
$$\bigwedge_{k = 1}^{K} \mathcal{B}_{n_k}\left( \pi_k  \right)  \underset{n_{K}}{\wedge} \mathcal{NB}(\alpha,p) = \mathcal{NB}\left(\alpha, \frac{p \prod_{k=1}^{K} \pi_k}{1-p\left(1-\prod_{k=1}^{K} \pi_k\right)}\right).$$

\noindent Since $M$ of those $\pi_k$ are beta distributed, we have by Fubini's theorem that, for $\boldsymbol{a} = (a_1, \dots, a_M)$ and $\boldsymbol{b} = (b_1, \dots, b_M)$,
$$\bigwedge_{k = 1}^K \mathcal{P}^{[c_k]}_{n_k}\left( \theta_k, \gamma_k \right) \underset{n_K}{\wedge} \mathcal{NB}(\alpha,p) = \mathcal{NB}\left(\alpha, \frac{p \gamma \pi}{1-p\left(1-\gamma \pi\right)}\right) \underset{\pi}{\wedge} \mathcal{PB}(\boldsymbol{a},\boldsymbol{b}),$$

\noindent where $\gamma = \prod_{k=1}^{K-M} \pi_k$, $\pi = \prod_{k=1}^M \pi_k$ and $\pi \sim \mathcal{PB}(\boldsymbol{a},\boldsymbol{b})$, the product beta distribution. 
Noticing that 
$$\frac{p \gamma \pi}{1-p\left(1-\gamma \pi\right)} = \frac{\frac{p \gamma}{1-p(1-\gamma)} \pi}{1-\frac{p \gamma}{1-p(1-\gamma)}\left(1-\pi\right)} \hspace{3 pt},$$
we have again that
\begin{align*}
    \mathcal{NB}\left(\alpha, \frac{p \gamma \pi}{1-p\left(1-\gamma \pi\right)}\right) \underset{\pi}{\wedge} \mathcal{PB}(\boldsymbol{a},\boldsymbol{b}) &= \left[ \mathcal{B}_{n_M}\left( \pi  \right) \underset{n_M}{\wedge} \mathcal{NB}\left(\alpha, \frac{p \gamma}{1- p(1-\gamma)}\right) \right] \underset{\pi}{\wedge} \mathcal{PB}(\boldsymbol{a},\boldsymbol{b})\\[5pt]
    &= \left[ \bigwedge_{k = 1}^{M} \mathcal{B}_{n_k}\left( \pi_k  \right) \underset{\pi}{\wedge} \mathcal{PB}(\boldsymbol{a},\boldsymbol{b}) \right] \underset{n_M}{\wedge} \mathcal{NB}\left(\alpha, \frac{p \gamma}{1- p(1-\gamma)}\right)\\[5pt]
    &= \bigwedge_{k = 1}^{M} \mathcal{BB}_{n_k}\left( a_k, b_k  \right) \underset{n_M}{\wedge} \mathcal{NB}\left(\alpha, \frac{p \gamma}{1- p(1-\gamma)}\right).
\end{align*}\qed

\noindent For the next result, we require the following lemma.
\begin{lemma}[page 98 in \cite{valiquette_thesis}]
    \label{lemma:binom_annex}
    For any $z \in \mathbb{R}$, $\pi \in (0,1)$ and $K \in \mathbb{N}_+$,
    $$\sum_{i=0}^n \binom{n}{i} \pi^i (1-\pi)^{n-i} \pFq{K+1}{K}{-i, \boldsymbol{a}}{\boldsymbol{a} + \boldsymbol{b}}{z} = \pFq{K+1}{K}{-n, \boldsymbol{a}}{\boldsymbol{a} + \boldsymbol{b}}{\pi z}.$$
\end{lemma}

\subsection*{Proof of Proposition \ref{prop:marginal_nb}}

First, let us calculate the probability generating function, denoted by $G(z)$, of
$\bigwedge_{k = 1}^K \mathcal{BB}_{n_k}\left(a_k, b_k  \right)$ using a proof by induction.
For $K = 1$, it is well known that $G(z) = \pFq{2}{1}{-n_1, a_1}{a_1 + b_1}{1-z}$.
Suppose for $K \geq 1$ that the generating function of this composition is given by 
\begin{equation}
\label{eqn:beta_binom_annex_2}
    G(z) = \pFq{K+1}{K}{-n_K, \boldsymbol{a}}{\boldsymbol{a} + \boldsymbol{b}}{1-z}.
\end{equation}
Let us prove the result for the composition of $K+1$ beta-binomial distributions.
We have by conditioning on the last $K$ terms that
$G(z) = \mathbb{E}\left[ \pFq{K+1}{K}{-n_K, \boldsymbol{a}_{-(K+1)}}{\boldsymbol{a}_{-(K+1)} + \boldsymbol{b}_{-(K+1)}}{1-z} \right]$, where the expectation is taken on the last $\mathcal{BB}_{n_{K+1}}(\alpha_{K+1}, \beta_{K+1})$.
Since the beta-binomial is a binomial mixture, the use of Lemma \ref{lemma:binom_annex} leads to
\begin{align*}
    G(z) &= \mathbb{E}\left[ \pFq{K+1}{K}{-n_{K+1}, \boldsymbol{a}_{-(K+1)}}{\boldsymbol{a}_{-(K+1)} + \boldsymbol{b}_{-(K+1)}}{(1-z)\pi} \right]\\
    &= \int_0^1 \frac{\pi^{a_{K+1} - 1} (1-\pi)^{b_{K+1} -1}}{\mathrm{B}(a_{K+1}, b_{K+1})} \pFq{K+1}{K}{-n_{K+1}, \boldsymbol{a}_{-(K+1)}}{\boldsymbol{a}_{-(K+1)} + \boldsymbol{b}_{-(K+1)}}{(1-z)\pi}d\pi\\
    &= \pFq{K+2}{K+1}{-n_{K+1}, \boldsymbol{a}}{\boldsymbol{a} + \boldsymbol{b}}{1-z},
\end{align*}
where, in the last equality, we used the integral representation of the generalized hypergeometric function \citep[e.g. 16.5.2 in ][]{NIST}. 
In order to find the probability generating function of the full distribution,
we only need to take the expectation of \eqref{eqn:beta_binom_annex_2} with respect to $n_K \sim \mathcal{NB}(\alpha,p)$.
Let $z \in (-p^{-1}, p^{-1})$, then 
\begin{align*}
    \pFq{K+1}{K}{-n_K, \boldsymbol{a}}{\boldsymbol{a} + \boldsymbol{b}}{1-z} &= \sum_{m=0}^n \binom{n}{m} (z-1)^m \frac{(\boldsymbol{a})_m}{(\boldsymbol{a}+\boldsymbol{b})_m}\\
    &\leq \sum_{m=0}^n \binom{n}{m} (z-1)^m = z^m,
\end{align*}
which implies that
$$G(z) = \sum_{m=0}^\infty \frac{(\alpha)_m}{m!} p^m (1-p)^{\alpha}  \pFq{K+1}{K}{-m, \boldsymbol{a}}{\boldsymbol{a} + \boldsymbol{b}}{1-z}\leq \left(\frac{1-p}{1-pz}\right)^\alpha.$$
Using a similar argument, we can prove the $n$-th derivative of $G(z)$ converges.
To obtain the p.m.f., it suffices to evaluate the term $G^{(n)}(0)/n!$ using equation 16.3.1 in \cite{NIST}.
For the case where $p \in (0, 1/2)$, we have by equation 1.22 in \cite{Norlund} that the infinite sum in the p.m.f. converges to $\pFq{K+1}{K}{\alpha+n, \boldsymbol{a} +n\mathbf{1}}{\boldsymbol{a} + \boldsymbol{b} + n\mathbf{1}}{\frac{p}{p-1}}$.\qed

\subsection*{Proof of Proposition \ref{prop:marginal_dirac}}

Using a similar argument as in Proposition \ref{prop:absorb}, it can be shown that we can interchange the order of composition, i.e. 
$$\bigwedge_{k = 1}^K \mathcal{P}^{[c_k]}_{n_k}\left( \theta_k, \tau_k \right) \underset{n_K}{\wedge} \mathcal{D}_m = \left[\bigwedge_{k = 1}^M \mathcal{BB}_{n_k}\left( a_k, b_k \right)\right] \underset{n_M}{\wedge} \mathcal{B}_m(\gamma).$$
Combining equation \eqref{eqn:beta_binom_annex_2} and Lemma \ref{lemma:binom_annex} above, the probability generating function of $X$ is given by 
$$G(z) = \mathrm{E}\left[ \mathrm{E}\left[ z^X \mid n_M \right] \right] = \mathrm{E}\left[ \pFq{M+1}{M}{-n_M, \boldsymbol{a}}{\boldsymbol{a} + \boldsymbol{\beta}}{1-z} \right] = \pFq{M+1}{M}{-m, \boldsymbol{a}}{\boldsymbol{a} + \boldsymbol{b}}{\gamma(1-z)} \hspace{3 pt}.$$
Since $\mathbb{P}(X=n) = G^{(n)}(0)/n!$, we can conclude.\qed

\subsection*{Proof of Proposition \ref{prop:factorial_mom}}

For each $\{j\} \in \mathfrak{L}$, there is a $\mathrm{Path}_{\{j\}}$ of length $K_{\{j\}}$. 
With each path, we can identify the leaves with the greatest path length.
Let us note that there are at least two leaves with maximum length due to the partition tree structure.
Furthermore, those leaves can be regrouped as siblings.
Without loss of generality, suppose the $m \geq 2$ first leaves have maximum length, and are siblings with a common parent node $A = \{1,\dots, m\}$.
From this set, we can use the law of total expectation yielding\vspace{-20 pt}

\begin{align*}
    \mathrm{E}[(\mathbf{Y})_{(\mathbf{r})} ] &= \mathrm{E}\left[(\mathbf{Y}_{-A})_{(\mathbf{r}_{-A})}  \mathrm{E}\left[ (\mathbf{Y}_{A})_{(\mathbf{r}_{A})}  
 \hspace{2pt} \bigg | \hspace{2pt} \mathbf{Y}_{-A}, \hspace{4pt} \abs{\mathbf{Y}_{A}} \right]\right].
\end{align*}
   
Since $\mathbf{Y}_A$ is conditionally independent of $\mathbf{Y}_{-A}$, and the distribution of $\mathbf{Y}_A$ given $\abs{\mathbf{Y}_{A}}$ is Pólya with parameter $\boldsymbol{\theta}_A$, then by Proposition \ref{prop:factorial_moment}   
$$\mathrm{E}[(\mathbf{Y})_{(\mathbf{r}})] = \frac{\prod_{C \in \mathfrak{C}_A}(\theta_C)_{(r_C, c_A)}}{(\abs{\boldsymbol{\theta}_A})_{(\abs{\mathbf{r}_A}, c_A)}} \mathrm{E}\left[ (\abs{\mathbf{Y}_A})_{(\abs{\mathbf{r}_A})} (\mathbf{Y}_{-A})_{(\mathbf{r}_{-A})} \right].$$

From this point on, the factorial moment of the full Tree Pólya can be obtained by calculating the factorial moment of a new Tree Pólya Splitting model where the $m$ first leaves are replaced by the leaf $A$, and calculating its $\abs{\mathbf{r}_A}$-th factorial moment. 
By iterating this process for the next set of siblings with maximal path length, we get the product over all internal nodes as mentioned.
Once this process arrives at the root $\Omega$, the factorial moment becomes
$$\mathrm{E}[(\mathbf{Y})_{(\mathbf{r})}] = \mathrm{E}[(\abs{\mathbf{Y}})_{(|\mathbf{r}|)}] \prod_{A \in \mathfrak{I}} \frac{\prod_{C \in \mathfrak{C}_A}(\theta_C)_{(\abs{\mathbf{r}_C}, c_A)}}{(\abs{\boldsymbol{\theta}_A})_{(\abs{\mathbf{r}_A}, c_A)}},$$
and the right-hand side expectation is simply $\mu_{\abs{\mathbf{r}}}$.\qed

\subsection*{Proof of Proposition \ref{prop:covariance_tree}}

By hypothesis, at least one path from a leaf to the ancestor node has a strictly positive length.
Otherwise, both $Y_i$ and $Y_j$ are siblings.
Without loss of generality, let us suppose that $\mathrm{Path}_{\{i\}}^{C_i}$ has length $K > 0$.
Then, for $A_k \in \mathrm{Path}_{\{i\}}^{C_i}$, we have by Corollary \ref{corollary:factorial_mom} that
$$\mathrm{E}[Y_i] = \left( \prod_{k=1}^{K} \frac{\theta_{A_{k-1}}}{\abs{\boldsymbol{\theta}_{A_k}}}\right) \mathrm{E}\big[\abs{\mathbf{Y}_{C_i}}\big].$$

Secondly, by a similar argument from the previous proof, we have
$$\mathrm{E}[Y_i Y_j] = \left( \prod_{k=1}^{K} \frac{\theta_{A_{k-1}}}{\abs{\boldsymbol{\theta}_{A_k}}}\right) \mathrm{E}\big[\abs{\mathbf{Y}_{C_i}} \hspace{1pt} Y_j \big].$$

By definition, the covariance is given by 
$$\mathrm{Cov}(Y_i, Y_j) = \left( \prod_{k=1}^{K} \frac{\theta_{A_{k-1}}}{\abs{\boldsymbol{\theta}_{A_k}}}\right) \mathrm{Cov}(\abs{\mathbf{Y}_{C_i}}, Y_j) = \frac{\gamma_i}{\gamma_{C_i}} \mathrm{Cov}(\abs{\mathbf{Y}_{C_i}}, Y_j)$$

By definition of $\gamma_\ell$.
Finally, if $\mathrm{Path}_{\{j\}}^{C_j}$ has length $0$, we conclude.
Otherwise, a similar argument on the expectations for $\mathrm{Path}_{\{j\}}^{C_j}$ yields the result.\qed

\subsection*{Proof of Proposition \ref{prop:covariance_tree_2}}

The covariance of $\abs{\mathbf{Y}_{C_i}}$ and $\abs{\mathbf{Y}_{C_j}}$ depends only on the parameter $\boldsymbol{\theta}_S$ and the first two factorial moments of $\abs{\mathbf{Y}_S}$, denoted for this proof by $\Tilde{\mu}_1$ and $\Tilde{\mu}_2$, at the ancestor node $S$.
Using Corollary \ref{corollary:factorial_mom} and the definitions of $\gamma_S$ and $\delta_S$ yields $\Tilde{\mu}_1 = \gamma_S \mu_1$, $\Tilde{\mu}_2 = \delta_S \gamma_S \mu_2.$
From equation \eqref{eqn:covariance}, the covariance of the subsums is given by $$\mathrm{Cov}(\abs{\mathbf{Y}_{C_i}}, \abs{\mathbf{Y}_{C_j}}) = \gamma_{C_i} \gamma_{C_j} \left[\frac{\abs{\boldsymbol{\theta}_S}}{(\abs{\boldsymbol{\theta}_S}+c)} \frac{\delta_S}{\gamma_S} \mu_{2} - \mu_{1}^2 \right].$$
We can conclude by combining this equality with Proposition \ref{prop:covariance_tree}.\qed

\subsection*{Proof of Proposition \ref{prop:correlation_tree}}

We only need to calculate the variance of $Y_\ell$ for $\ell = i,j$. 
Again, by Corollary \ref{corollary:factorial_mom}, the result follows since
$$\mathrm{Var}(Y_\ell) = \mathrm{E}[Y_\ell(Y_\ell-1)] + \mathrm{E}[Y_\ell]\big(1-\mathrm{E}[Y_\ell]\big) = \gamma_\ell \big( \delta_\ell \mu_2 + \mu_1 (1-\gamma_\ell \mu_1) \big).$$\qed

\bibliography{biblio}

\end{document}